\documentclass[a4paper,11pt,english]{amsart}
\usepackage{amsmath,amsthm,amsfonts,amssymb,upgreek,multicol,mathrsfs,pifont,amscd,latexsym,cite,bm,geometry,color,fancyhdr,abstract,framed,appendix}
\usepackage{graphicx,algorithm,algorithmic,tabularx,url,float,booktabs,array,threeparttable,tikz}

\newcommand{\mn}{\mathbf{n}}
\newcommand{\mt}{\mathcal{T}_h}

\newcommand{\me}{\mathcal{E}}
\newcommand{\ms}{\mathbb{S}}

\newcommand{\mbb}{\mathbb{R}^{d}}
\newcommand{\be}{\begin{eqnarray}}
\newcommand{\ee}{\end{eqnarray}}
\newcommand{\ea}{\end{array}}
\newcommand{\ba}{\begin{array}}
\newcommand{\bsi}{\bm{\sigma}}
\newcommand{\bta}{\bm{\tau}}
\newcommand{\rd}{\rm{div}}
\newcommand{\bd}{\mathbf{div}}
\newcommand{\bc}{\mathbf{curl}}

\newcommand{\ben}{\begin{equation*}}
\newcommand{\een}{\end{equation*}}

\newtheorem{theorem}{Theorem}[section]

\newtheorem{lemma}{Lemma}[section]

\newtheorem{remark}{Remark}[section]

\numberwithin{equation}{section}
\numberwithin{table}{section}
\numberwithin{figure}{section}

\begin{document}



\title[]{
A family of mixed finite elements for the biharmonic equations on triangular and tetrahedral grids
}

\author {Jun Hu}
\address{LMAM and School of Mathematical Sciences, Peking University,
  Beijing 100871, P. R. China. }
\email{  hujun@math.pku.edu.cn}

\author {Rui Ma}
\address{Universit\"at Duisburg-Essen, Thea-Leymann-Str. 9, {\rm{45127}} Essen, Germany. }
\email{ rui.ma@uni-due.de }

 \author {Min Zhang}
\address{School of Mathematical Sciences, Peking University,
	Beijing {\rm100871}, P. R. China.}
\address{Computational Science Research Center, Beijing {\rm{100193}}, P. R. China.}
\email{  zmzoe@pku.edu.cn }

\thanks{The first author was supported by the NSFC Projects 11625101 and 11421101.}
\maketitle

\begin{abstract}
This paper introduces a new family of mixed finite elements for solving a mixed formulation of the biharmonic equations in two and three dimensions. The symmetric stress $\bsi=-\nabla^{2}u$ is sought in the Sobolev space $H({\rd}\bd,\Omega;\ms)$ simultaneously with the displacement $u$ in $L^{2}(\Omega)$. Stemming from the structure of $H(\bd,\Omega;\ms)$ conforming elements for the linear elasticity problems proposed by J. Hu and S. Zhang, the $H({\rd}\bd,\Omega;\ms)$ conforming finite element spaces are constructed by imposing the normal continuity of $\bd\bsi$ on the $H(\bd,\Omega;\ms)$ conforming spaces of $P_{k}$ symmetric tensors. The inheritance makes the basis functions easy to compute. The discrete spaces for $u$ are composed of the piecewise $P_{k-2}$ polynomials without requiring any continuity. Such mixed finite elements are inf-sup stable on both triangular and tetrahedral grids for $k\geq 3$, and the optimal order of convergence is achieved. Besides, the superconvergence and the postprocessing results are displayed. Some numerical experiments are provided to demonstrate the theoretical analysis.
\end{abstract}



\maketitle
\def\a#1{\begin{align*}#1\end{align*}}\def\an#1{\begin{align}#1\end{align}}
\def\ad#1{\begin{aligned}#1\end{aligned}}

\section{Introduction}
Let $\Omega\subset\mbb$ be a bounded Lipschitz polyhedral domain with $d=2$ or $3$. Given a load $f\in L^2(\Omega)$, consider the biharmonic equation 
\be
\label{biha}
\left\{
\begin{array}{rll}
\Delta^2 u=f,  &in~ \Omega,\\
u=u_{n}=0, &on ~ \partial\Omega.\\
\end{array}\right.
\ee
Here $\Delta^2$ is the biharmonic operator, $\mn$ is the unit outer normal to the boundary $\partial \Omega$, and $u_{n}:=\partial u/\partial \mn$. 

Many attempts have been made to approach the biharmonic problem \eqref{biha}, ranging from conforming and classical
nonconforming finite element methods, discontinuous Galerkin methods to mixed methods, such as \cite{Susanne2010A, HuSpp,  Lascaux1975Some, Ciarlet1974A, Gerasimov2012CORNERS, Behrens2011A, Herrmann, Johnson1973On}, to name just a few. 
On triangular grids,  the lowest order of polynomials of the  $H^{2}$ conforming finite elements is 5. That is the Argyris element \cite{Arg,Ciarletbook}, and it can be reduced to the Bell element \cite{wangshi,Ciarletbook} with $18$ degrees of freedom.
On tetrahedral grids, a $P_{9}$ element constructed in \cite{Zhang2009A} is the lowest order conforming element. In general, due to the high degrees of freedom with higher order derivatives of the $H^{2}$ conforming elements, in addition to the complexity in construction, the computation is relatively costing. Nevertheless, some conforming finite elements are developed \cite{Ciarletbook,wangshi,DDPS,PS,FVS,ZA,2011The, 2015The}.
One way to reduce the high degrees of freedom is to use nonconforming finite elements, such as the Morley element \cite{Morley(1968), Ciarletbook,wangshi}, the Adini element \cite{AdiniClough,Ciarletbook,wangshi},  the Veubake element \cite{VB}, a class of Zienkiewicz-type nonconforming elements in any dimensions designed in \cite{NTZ}, and other higher order nonconforming methods \cite{BSM,WZZ,CCZ,GJN,HuSpp,article}. The other way is to adopt different variational principles to avoid computational difficulty. A popular choice is mixed finite element methods. For example, the Ciarlet-Raviart method \cite{Ciarlet1974A} turns \eqref{biha} into a lower order system by introducing an auxiliary variable $\phi=-\Delta u$, and casts the new system in variational form, then considers the Ritz-Galerkin method corresponding to this variational formulation. {However, such decoupling may not be valid if the polygonal domain is not convex; see \cite{MR2479119}}. Instead of $\phi=-\Delta u$,  the matrix of the second partial derivatives of $u$, ${\bsi}=-\nabla^{2}u$ is introduced in the Hermann-Miyoshi method \cite{Herrmann, Miyoshi}. A further mixed method for \eqref{biha} is the Hermann-Johnson element, and the auxiliary variable introduced is the same as the Hermann-Miyoshi method, while the continuity of $\mathbf{n}^{{\rm{T}}}\bsi\mathbf{n}$ is imposed on $\bsi$.

In this paper, a more intrinsic variational formulation is considered, {and it is also known as the Hodge-Laplacian boundary value problem of the divdiv complex. In \cite{2020Complexes}, the well-posedness of the Hodge-Laplacian boundary value problem is discussed. The mixed finite element method seeks the stress ${\bsi}=-\nabla^{2}u$ in the }Sobolev space $H({\rd}\bd, \Omega; \ms)$ with 
\be
H({\rd}\, {\bd}, \Omega; \ms):=\{\bta \in L ^{2}(\Omega;\ms):\, {\rd}\, {\bd}\bta\in L ^{2}(\Omega)\},
\ee
equipped with the squared norm
\be
\label{divNorm}
\|\bta\|_{H({\rd}{\bd})}^{2}:=\|\bta\|_{0}^{2}+\|{\rd}{\bd}\bta\|_{0}^{2}.
\ee
Here $\ms$ denotes the set of symmetric $\mathbb{R}^{d\times d}$ matrices.
{Simultaneously, the mixed method seeks} $u\in L^{2}(\Omega)$ such that
\be
\begin{aligned}
\label{ConV}
(\bsi,\bta)+({\rd} {\bd}\, \bta, u) =& 0 & \text{for all}\,  \bta\in H(\rm{div}\, \mathbf{div}, {\rm{\Omega}}; \ms),\\
({\rd}\, \bd\bsi, v) = &-(f, v) &\text{for all} \, v\in L^{2}(\Omega).
\end{aligned}
\ee
It is not easy to construct a $H({\rd}\, \bd, {\rm{\Omega}}; \ms)$ conforming element, and the symmetry of the tensor makes things more complex. 
A family of $H(\bd;\ms)$ conforming finite elements for elasticity equations is proposed in \cite{Hu2014A, Hu2015FINITE, MR3301063}. If $\bd\bsi\in H({\rd})$ holds for all $\bsi\in H(\bd;\ms)$, then $\bsi\in H({\rd}\bd;\ms)$ follows. The relation triggers an idea to obtain the $H({\rd}\, \bd, {\rm{\Omega}}; \ms)$ conforming elements by imposing the continuity of $\mathbf{n}^{{\rm{T}}}\bd\bsi$ on $H(\bd;\ms)$ conforming spaces. A question arises naturally how to characterize this additional continuity appropriately.  

Attempts have been made in \cite{Yxq2,Yxq1}, where the stress space is composed by the aforementioned $H(\bd;\ms)$ conforming elements \cite{Hu2014A, Hu2015FINITE, MR3301063}, and the displacement space chooses the $P_{k}$ conforming finite element with $k \geq 2$. However, the $L^{2}$ norms are not optimal. {In \cite{MR3864098,MR3925478}, a depiction of the Sobolev space $H({\rd}\, \bd, \Omega; \ms)$ is introduced, and the discontinuous Petrov-Galerkin method is considered. }
Recently, some finite element spaces for $H({\rd}\, \bd, \Omega; \ms)$ conforming symmetric tensors are constructed on triangles \cite{Chen2020Finite} and tetrahedrons \cite{Chen2020Finite3d}. These elements are exploited to solve the mixed problem \eqref{ConV} and the optimal order of convergence is achieved.
In two dimensions, a simple application of Green's formula shows
\ben
\begin{split}
~~~~({\rd}\bd\bsi, v)_{K}&=(\bsi,\nabla^{2} v)_{K} +\sum_{e\in \mathcal{E}(K)}(\mathbf{n}^{{\rm{T}}}\bd\bsi,v)_{e}-\sum_{e\in\mathcal{E}(K)}(\bsi\mathbf{n}, \nabla v)_{e}, \\
\end{split}
\een
{with the unit out normal vector $\mn=(n_{1},n_{2})^{{\rm{T}}}$ and the unit tangent vector $\mathbf{t}=(-n_{2},n_{1})^{{\rm{T}}}$ below. }
Expand  $(\bsi\mathbf{n}, \nabla v)_{e}=(\mathbf{n}^{{\rm{T}}}\bsi\mathbf{n},\partial_{n} v)_{e}+(\mathbf{t}^{{\rm{T}}}\bsi\mathbf{n},\partial_{t} v)_{e}$. A further integration by parts gives rise to 
\be
\label{green:1}
\begin{split}
~~~~({\rd}\bd\bsi, v)_{K}&=(\bsi,\nabla^{2} v)_{K} -\sum_{e\in \mathcal{E}(K)}\sum_{a\in \partial{e}}{\rm{sign}}_{e,a}(\mathbf{t}^{{\rm{T}}}\bsi\mathbf{n})(a)v(a)\\
&-\sum_{e\in\mathcal{E}(K)}[(\mathbf{n}^{{\rm{T}}}\bsi\mathbf{n}, \partial_{n} v)_{e}-(\partial_{t}(\mathbf{t}^{{\rm{T}}}\bsi\mathbf{n})+\mathbf{n}^{{\rm{T}}}\bd\bsi,v)_{e}],\\
\end{split}
\ee
with
\ben
{\rm{sign}}_{e,a}:=\left\{\begin{array}{ll}
1,~~~~&\text{if {$a$} is the end point of $e$},\\
-1,~~~~&\text{if  {$a$} is the start point of $e$}.\\
\end{array}\right.
\een 
Based on \eqref{green:1}, besides the normal-normal continuity, the stress tensor is continuous at vertices and another trace involving the combination of derivatives of the stress is identified.
The basic design of the $H({\rd}\bd, \Omega; \ms)$ conforming finite elements in \cite{Chen2020Finite} follows.

However, it is arduous to compute the basis functions for the elements in \cite{Chen2020Finite, Chen2020Finite3d}.  Motivated by \cite{Hu2014A, Hu2015FINITE,MR3301063}, {this} paper introduces a more straight forward characterization of the {$H({\rd}\bd;\ms)\cap H({\bd};\ms)$} space. Instead of involving combination of derivatives of stresses, the continuity of $\bsi\mathbf{n}$ and $\mathbf{n}^{{\rm{T}}}\bd\bsi$ is imposed in the design of the new $H({\rd}\bd;\ms)$ conforming elements. {Correspondingly, the finite elements obtained in this paper are more regular than those in \cite{Chen2020Finite, Chen2020Finite3d}. Actually, {the new elements} are subspaces} {of the elements proposed in \cite{Chen2020Finite, Chen2020Finite3d}}. 
The $H(\bd;\ms)$ bubble functions presented in \cite{Hu2014A, Hu2015FINITE, MR3301063} possess vanishing $\bsi\mathbf{n}$ on each face. Therefore, the basis functions corresponding to the degrees of freedom $\mathbf{n}^{{\rm{T}}}\bd\bsi$ can be expressed linearly by the basis of these bubbles. The remainder basis functions can be derived by the former $\mathbf{n}^{{\rm{T}}}\bd\bsi$ basis and the basis functions given by \cite{Hu2014A, Hu2015FINITE, MR3301063}. Besides, the new $H({\rd}\bd;\ms)$ conforming finite elements in two and three dimensions can be constructed in an almost unified way, while the degrees of freedom in \cite{Chen2020Finite3d} are fairly sophisticated. 

{In addition, a vectorial $H^{1}$ conforming finite element on triangular grids is introduced, and this element plus the $H({\rd}{\bd};\ms)$ conforming finite element form the discrete divdiv complex. In this paper, the exactness of the finite element analogy of divdiv complex is proved on a contractible domain. Actually, by rotation, the {two dimensional} divdiv complex is equivalent to the strain complex. Conforming finite elements for $H({\rm{rot}\mathbf{rot}};\ms)$ are obtained in \cite{Chen2020Finite} in two dimensions. By using piecewise polynomials based on the Clough-Tocher split of the triangle, some lower-order $H({\rm{rot}\mathbf{rot}};\ms)$ conforming finite elements are constructed to obtain the discrete strain complex in \cite{2019FiniteHu}.}

Furthermore, the new  $H({\rd}\, {\bd}, \Omega; \ms)$ conforming finite elements space developed for $d$ being $2$ and $3$ are capable of discretizing the mixed formulation \eqref{ConV} with the optimal order of convergence.  

The remainder of this paper is organized as follows. In the subsequent section, the construction of $H({\rd}\bd;\ms)$ conforming finite elements in two dimensions as well as in three dimensions is presented. {Correspondingly, a vectorial $H^{1}$ conforming finite element in two dimensions is introduced to establish the discrete complex, which is proved to be exact on a contractible domain}.  In Section 3, the new conforming elements are exploited to discrete the mixed problem \eqref{ConV}. The well-posedness is proved and the error analysis follows. Besides, superconvergence and postprocessing results are displayed. In Section 4,  numerical examples are presented to demonstrate the theoretical analysis results. In the end, the appendix provides some ideas to construct the basis functions by a specific example. 

Throughout the paper, an inequality $\alpha\lesssim \beta$ replaces $\alpha\leq c  \beta$ with some multiplicative mesh-size independent constant $c>0$, which depends on $\Omega$ only. While $\alpha\thicksim  \beta$ means $\alpha\lesssim  \beta$ and $ \beta\lesssim \alpha$ hold simultaneously. Standard notation on Lebesgue and Sobolev spaces are employed. For a subset $G\subset \Omega$, $(\cdot,\cdot)_{G}$ denotes the $L^2$ scalar product over $G$, $\Vert{\cdot}\Vert_{0,G}$ denotes the $L^2$ norm over a set $G$.  $\Vert{\cdot}\Vert_0$ abbreviates $\Vert{\cdot}\Vert_{0,\Omega}$. Other cases are similar. Let $\mathcal{D}(G)$ denote the set of all infinitely differentiable compactly supported functions on $G$. {Let $P_l(G)$ stand for the set of all polynomials with the total degree no more than $l$ over $G$. 
Notation $\mathbb{X}$ could be $\mathbb{R}$, $\mathbb{R}^{d}$, $\mathbb{M}$, $\mathbb{T}$, $\ms$, and $\mathbb{K}$ in the text. Correspondingly, it denotes the space of scalars, vectors in $d$ dimensions, matrices in $\mathbb{R}^{d\times d}$, traceless matrices in $\mathbb{R}^{d\times d}$, symmetric matrices in $\mathbb{R}^{d\times d}$, and skew-symmetric matrices in $\mathbb{R}^{d\times d}$, respectively. Dimension $d$ is either $2$ or $3$ in this paper, and it coincides with shape of $G$. For instance, 
any variable in $\mathcal{D}(G;\ms)$ is a symmetric matrix on $G$ and it is infinitely differentiable compactly supported. } Similarly, {$P_{l}(G;\mathbb{X})$ can be defined in the same way}. Generally, $\mathcal{D}(G;\mathbb{R})$ is simply abbreviated as $\mathcal{D}(G)$, and so does $P_{l}(G)$ for $P_{l}(G;\mathbb{R})$. 
Denote the $\bc$ operators below,
 \ben
 \begin{aligned}
 &\bc \varphi=(-\partial_{y}\varphi, \partial_{x}\varphi)^{{\rm{T}}}&\text{for all} \,  \varphi\in \mathcal{D}(\Omega;\mathbb{R}),\, d=2.\\
&\bc \bm{\varphi}=(\partial_{y}\varphi_{3}-\partial_{z}\varphi_{2}, \partial_{z}\varphi_{1}-\partial_{x}\varphi_{3}, \partial_{x}\varphi_{2}-\partial_{y}\varphi_{1})^{{\rm{T}}}&\text{for all}\, \bm{\varphi}\in \mathcal{D}(\Omega;\mathbb{R}^{3}),\, d=3.\\
 \end{aligned}
 \een
{The symmetric gradient operator is denoted as $\bm{\varepsilon}(u)=1/2(\nabla u+(\nabla u)^{{\rm{T}}})$}. Generally, for a column vector function, differential operators for scalar functions will be applied row-wise to produce a matrix function. Similarly for a matrix function, differential operators for vector functions are applied row-wise. However $\bc^{*}$ will be the $\bc$ operator applied column-rise.

\section{The conforming finite element spaces}
This section covers some preliminaries and the construction of the new $H({\rd}\, {\bd}, \Omega; \ms)$ conforming finite elements in both two and three dimensions. Besides, a vectorial $H^{1}(\Omega;\mathbb{R}^{2})$ conforming finite element space is introduced, and a discrete case of Hilbert complex is obtained.
\subsection{Notation}Suppose $\mathcal{T}_h$ is a shape regular subdivision of $\Omega$ consisting of triangles in two dimensions and tetrahedrons in three dimensions. Denote $h$ the maximum of the diameters of all elements $K\in \mt$. Let $\mathcal{E}_{h}$, $\mathcal{F}_{h}$ and $\mathcal{V}_{h}$ be the set of all edges, faces, and vertices of $\Omega$ regarding to $\mt$, respectively. Given $K\in\mt$, let $\mathcal{E}(K)$ denote the set of all edges of $K$, and $h_{e}$ stands for the diameter of edge $e\in \mathcal{E}_{h}$. Furthermore, when $d=3$, define the set of all facets of the tetrahedron $K$ as $\mathcal{F}(K)$, and $h_{F}$ stands for the diameter of face $F\in \mathcal{F}_{h}$.
Let $\mathbf{n}$ and $\mathbf{t}$ be the unit outer normal and unit tangential vector of $\partial K$ respectively. More specific, when $d=2$, ${\mathbf{t}_{e}}$ denotes the unit tangential vector along $e\in\mathcal{E}(K)$, and ${\mathbf{n}_{e}}$ is the normal counterpart. While $d=3$, given $e\in \mathcal{E}(K)$, the unit tangential vector $\mathbf{t}_{e}$, as well as two unit normal vectors, ${\mathbf{n}_{e, 1}}$ and ${\mathbf{n}_{e, 2}}$ are fixed.  For a facet $F\in \mathcal{F}(K)$, the unit outer normal vector ${\mathbf{n}_{F}}$ as well as two unit tangential vectors ${\mathbf{t}_{F, 1}}$ and ${\mathbf{t}_{F, 2}}$ are fixed. Within the context, $\mathbf{t}_{i}$ and $\mathbf{n}_{i}$  abbreviate $\mathbf{t}_{F, i}$ and $\mathbf{n}_{e, i}$, respectively, $i=1,2$. Besides, the union of all vertices of $K$ is denoted as $\mathcal{V}(K)$.
The jump of $u$ across an interior $d-1$ face $G$ shared by neighboring elements $K_+$ and $K_-$ is defined by
$$\left [u\right]_G:=u|_{K_+}-u|_{K_-}.$$
When it comes to any boundary face $G\subset\partial \Omega$, the jump $[\cdot]_G$ reduces to the trace.  
 
 For ensuing analysis, let $RM(K)$ denote local rigid motions. When $K$ is a triangle with ${\mathbf{x}}=(x, y)^{{\rm{T}}}\in {K}$, 
 \be
 RM_{\triangle_{2}}(K)=\left\{\begin{pmatrix}
 c_{1}+c_{3}y\\
 c_{2}-c_{3}x\\
 \end{pmatrix} : c_{1}, c_{2}, c_{3}\in \mathbb{R}
  \right\}.
 \ee
If $K$ is a tetrahedron with $\bf{x}=(\emph{x}, \emph{y},\emph{z})^{{\rm{T}}}\in \emph{K}$, then
 \be
 RM_{\triangle_{3}}(K)=\left\{\begin{pmatrix}
 c_{1}-c_{4}y-c_{5}z\\
 c_{2}+c_{4}x-c_{6}z\\
 c_{3}+c_{5}x+c_{6}y\\
 \end{pmatrix} : c_{1}, c_{2}, c_{3}, c_{4}, c_{5}, c_{6} \in \mathbb{R}
  \right\}.
 \ee
 
Define
\ben
P_{h}:=\{q\in L^{2}(\Omega):q|_{K}\in P_{k-2}(K)~\text{for all\,} K\in \mt\}.
\een
Denote $P_{h}$ as $P_{h,\triangle_{2}}$ or $P_{h,\triangle_{3}}$ in two and three dimensions respectively.

Besides, $RT_{k}$ is the Raviart-Thomas element space \cite{brezzi2012mixed},
\ben
RT_{k}(K;\mathbb{R}^{d})=P_{k}(K; \mathbb{R}^{d})+\mathbf{x}P_{k}(K).
\een
Notice that
\ben
{\rm{dim}} \, RT_{k}(K;\mathbb{R}^{3})=\frac{(k+1)(k+2)(k+4)}{2}.
\een
Denote $RT$ the lowest order  Raviart-Thomas element space on $\Omega$.

\subsection{The construction of the conforming elements on triangular grids}On each triangle $K$, denote  $\lambda_{i}$, $i=1,2,3$ the barycenter coordinates. The finite element shape functions are simply formed by $P_{k}(K; \ms)$, $k\geq 3$.  Some results  are presented in the following two lemmas for later use. 
\begin{lemma}[\cite{MfemB}]
\label{2d:5}
Given $K\in \mt$, suppose ${\bm{\psi}}\in P_{k}(K;\mathbb{R}^{2})$ satisfies ${\rd}{\bm{\psi}}=0$, and ${\bm{\psi}}\cdot\mathbf{n}\, |_{\partial K}=0$. Then there exists some $q\in  \lambda_{1} \lambda_{2} \lambda_{3}P_{k-2}(K)$, such that $${\bm\psi}={\mathbf{curl}} \, q.$$ 
\end{lemma}

\begin{lemma}[\cite{Winther2002Mixed,2018Nodal}]
\label{2d:6}
Given  $K\in \mt$, suppose $\bta\in P_{k}(K;\ms)$ satisfies $\bd\bta=0$, and $\bta\mathbf{n}|_{\partial K}=0$. Then there exists some $q\in  (\lambda_{1}\lambda_{2}\lambda_{3})^{2}P_{k-4}(K)$, such that $$\bta=\mathcal{J}q$$ with
\be
\mathcal{J} q:=
\begin{pmatrix}
\frac{\partial^{2}q}{\partial y^{2}} &-\frac{\partial^{2}q}{\partial x\, \partial y}\\
-\frac{\partial^{2}q}{\partial x\, \partial y} &\frac{\partial^{2}q}{\partial x^{2}}\\
\end{pmatrix}.
\ee
\end{lemma}

The degrees of freedom are defined as follows.
\begin{align}
   \bsi (a)&~~ \text{for all}\,  a\in \mathcal{V}(K); \label{2dDof:1}\\
 (\bsi \mathbf{n}, \bm{\phi})_{e}& ~~\text{for all}\, \bm{\phi} \in P_{{k-2}}(e; \mathbb{R}^{2}), \, e\in \mathcal{E}(K); \label{2dDof:2}\\
 (\bd \bsi\cdot \mathbf{n},\,  q)_{e}&~~ \text{for all}\, q \in P_{{k-1}}(e), \, e\in \mathcal{E}(K); \label{2dDof:3}\\
 (\bsi, \nabla^{2} \,q)_{K}& ~~ \text{for all}\, q \in P_{{k-2}}(K); \label{2dDof:4}\\
 (\bsi, \nabla \bc \, q)_{K}& ~~\text{for all}\, q \in \lambda_{1}\lambda_{2}\lambda_{3} P_{k-3}(K)/P_{0}(K);  \label{2dDof:5}\\
  (\bsi, \mathcal{J} q)_{K}&~~\text{for all}\, q\in  (\lambda_{1}\lambda_{2}\lambda_{3})^{2}P_{k-4}(K).\label{2dDof:6}
\end{align}

{
\begin{remark}
Any function $\varphi \in \lambda_{1}\lambda_{2}\lambda_{3} P_{k-3}(K)/P_{0}(K)$ means $\varphi= \lambda_{1}\lambda_{2}\lambda_{3} q$ for some $q\in P_{k-3}(K)$ as well as $\int_{K}\varphi=0$.
\end{remark}}

Degrees of freedom \eqref{2dDof:1}--\eqref{2dDof:3}  characterize the continuity of the space $H({\rd}\bd;\ms)$. 
With the help of Lemma \ref{2d:5} and  Lemma \ref{2d:6}, \eqref{2dDof:5}--\eqref{2dDof:6} can be used to derive the unisolvence. Besides, the degrees of freedom \eqref{2dDof:1}--\eqref{2dDof:2} are exactly the characterization of the continuity of $H({\bd}; \ms)$ in \cite{Hu2014A, Hu2015FINITE}, and the continuity of \eqref{2dDof:3} across edges leads to $\bd\bsi\in H({\rd}; \mathbb{R}^{2})$.

The global finite element space is defined by
\be
\label{2d:space}
\begin{split}
\Sigma_{k,\triangle_{2}}:=&\{\bm{\tau}\in H({\rd}{\bd},\Omega;\ms):\, \bm{\tau}|_{K}\in P_{k}(K;\ms) ~~\text{for all}~ K\in\mt, \\
& \text{all the degrees of freedom \eqref{2dDof:1}--\eqref{2dDof:6} are single-valued}\}.
\end{split}
\ee
\begin{theorem}
The degrees of freedom \eqref{2dDof:1}--\eqref{2dDof:6} uniquely determine a polynomial of $P_{k}(K; \ms)$ in the space $\Sigma_{k,\triangle_{2}}$ defined in \eqref{2d:space}.
\end{theorem}
\begin{proof}
To start with, it is easy to check that the number of the degrees of freedom \eqref{2dDof:1}--\eqref{2dDof:6} is {
\ben
\begin{split}
&9+6(k-1)+3k+\frac{k(k-1)}{2}-3+\frac{(k-1)(k-2)}{2}-1+\frac{(k-2)(k-3)}{2}\\
&=\frac{3(k+1)(k+2)}{2}={\rm{dim}}P_{k}(K; \ms).
\end{split}
\een
} 
It suffices to prove if degrees of freedom \eqref{2dDof:1}--\eqref{2dDof:6} vanish for $\bsi\in P_{k}(K;\ms)$, then $\bsi=0$.
Given any $v\in P_{k-2}(K)$, integration by parts and the zero degrees of freedom \eqref{2dDof:2}--\eqref{2dDof:4} lead to 
\be
\label{Green}
\begin{split}
({\rd}\bd \bsi, v)_{K}&=(\bsi, \nabla^{2}\emph{v})_{K}-\sum_{e\in\mathcal{E}(K)}(\bsi \mathbf{n}, \nabla \emph{v})_{e}+\sum_{e\in\mathcal{E}(K)}({\rd} \bsi\cdot \mathbf{n}, v)_{e}\\
&=0.
\end{split}
\ee
This implies $\rd \bd\bsi=0$. 
Together with \eqref{2dDof:3}, according to Lemma \ref{2d:5},  there exists some $\varphi\in\lambda_{1}\lambda_{2}\lambda_{3} P_{k-3}(K)$ such that
\be
\label{teq}
\bd \bsi=\bc \varphi.
\ee
For any function $\vartheta\in \lambda_{1}\lambda_{2}\lambda_{3} P_{k-3}(K)/P_{0}(K)$, integration by parts plus \eqref{2dDof:2} and \eqref{2dDof:5} show
\be
\label{t1}
\begin{split}
(\bc\varphi,\bc \vartheta )_{K}&=(\bd\bsi, \bc\vartheta )_{K}\\
&=-(\bsi, \nabla \bc\vartheta)_{K}+\sum_{e\in\mathcal{E}(K)}(\bsi \mathbf{n}, \bc\vartheta)_{e}=0.
\end{split}
\ee
Besides, \eqref{2dDof:1}--\eqref{2dDof:2} result in 
\be
\label{trm}
(\bd \bsi, \mathbf{v})_{K}=0 ~~\text{for all}~ \mathbf{v}\in RM_{\triangle_{2}}(K).
\ee
Take $\mathbf{v}=(-y,x)^{{\rm{T}}}$ in \eqref{trm}, using \eqref{teq} for replacement, 
\be
\label{t2}
(\bd \bsi, \mathbf{v})_{K}=(\bc\varphi, \mathbf{v})_{K}=2\int_{K}\varphi \, d\mathbf{x}=0.
\ee
Note that \eqref{t1} and \eqref{t2} lead to $\varphi=0$, thus $\bd\bsi=\mathbf{0}$. Furthermore, due to \eqref{2dDof:1}--\eqref{2dDof:2}, according to Lemma \ref{2d:6}, $\bd\bsi=0$ entails the relation $\bsi=\mathcal{J}\zeta$ for some $\zeta\in (\lambda_{1}\lambda_{2}\lambda_{3})^{2}P_{k-4}(K)$. This and \eqref{2dDof:6} conclude $\bsi=0$ immediately. 
\end{proof}
\begin{remark}
The degrees of freedom $\bd\bsi\cdot\mathbf{n}$ in \eqref{2dDof:3} can be replaced by $\partial_{n}(\mathbf{n}^{{\rm{T}}}\bsi\mathbf{n})$ since
\ben
\bd\bsi\cdot \mathbf{n}={\bd(\bsi\mathbf{n})}=\partial_{t}(\mathbf{t}^{{\rm{T}}}\bsi\mathbf{n})+\partial_{n}(\mathbf{n}^{{\rm{T}}}\bsi\mathbf{n}).
\een
Let $a_{e,1}$ and $a_{e,2}$ be the start and end point of the edge $e$ respectively. Integration by parts leads to 
\ben
(\partial_{t}(\mathbf{t}^{{\rm{T}}}\bsi\mathbf{n}),v)_{e}=\mathbf{t}^{{\rm{T}}}\bsi\mathbf{n}\,v\big{|}_{a_{e,1}}^{a_{e,2}}-(\mathbf{t}^{{\rm{T}}}\bsi\mathbf{n},\partial_{t} v)_{e}~~\text{for all}~~ v\in P_{k-1}(e).
\een
The first term can be covered by the degrees of freedom \eqref{2dDof:1}, and the second term can be derived by the degrees of freedom \eqref{2dDof:2}.
\end{remark}

\subsection{The construction of {the finite element divdiv complex} on triangular grids}
Define
\ben
H^{1}({\rd},\Omega;\mathbb{R}^{2}):=\{\bm{v}\in H^{1}(\Omega;\mathbb{R}^{2}): {\rd}\bm{v}\in H^{1}(\Omega) \}.
\een

The vectorial space $V_{h}\subset {H^{1}({\rd}, \Omega;\mathbb{R}^{2})}$ is introduced in this subsection, and the discrete exact complex is established. On a triangle $K\in\mt$, the shape function space is $P_{k+1}(K;\mathbb{R}^{2})$, and the degrees of freedom are
\begin{align}
\bm{v}(a),{\nabla}{\bm{v}}(a)&~~  \text{for all}\, a\in \mathcal{V}(K); \label{H1Dof:1}\\
 ({\bm{v}}, {\bm{\phi}})_{e}&~~ \text{for all}\, {\bm{\phi}} \in P_{{k-3}}(e;\mathbb{R}^{2}), \, e\in \mathcal{E}(K); \label{H1Dof:2}\\
 ({\rd}{\bm{v}},  q)_{e}&~~ \text{for all}\, q \in P_{{k-2}}(e), \, e\in \mathcal{E}(K); \label{H1Dof:3}\\
 ({\bm{v}}, \nabla q)_{K}&~~ \text{for all}\, q \in P_{{k-3}}(K); \label{H1Dof:4}\\
 ({\bm{v}}, {\bc}\,  q)_{K}&~~ \text{for all}\, q \in (\lambda_{1}\lambda_{2}\lambda_{3})^{2}P_{k-4}(K). \label{H1Dof:5}
\end{align}

Then the space $V_{h}$ is defined by
\be
\label{H1:space}
\begin{split}
V_{h}:=&\{\bm{v}\in H^{1}(\Omega;\mathbb{R}^{2}):\, \bm{v}|_{K}\in P_{k+1}(K;\mathbb{R}^{2})~~\text{for all}~ K\in\mt, \\
& \text{all the degrees of freedom \eqref{H1Dof:1}--\eqref{H1Dof:5} are single-valued}\}.
\end{split}
\ee
\begin{theorem}
The degrees of freedom \eqref{H1Dof:1}--\eqref{H1Dof:5} uniquely determine a polynomial of $P_{k+1}(K; \mathbb{R}^{2})$ in the space $V_{h}$ defined in \eqref{H1:space}.
\end{theorem}
\begin{proof}
To start with, it is easy to check that the number of the degrees of freedom \eqref{H1Dof:1}--\eqref{H1Dof:5} equals to the dimension of $P_{k+1}(K; \mathbb{R}^{2})$. In fact, both of them are 
\ben
{(k+3)(k+2)}.\\
\een

It suffices to prove if  degrees of freedom \eqref{H1Dof:1}--\eqref{H1Dof:5} vanish for ${\bm{v}} \in P_{k+1}(K; \mathbb{R}^{2})$, then $\bm{v}=0$. 
Actually, \eqref{H1Dof:1}--\eqref{H1Dof:2} lead to 
\be
\label{H1:2a}
\bm{v}|_{e}=0 ~~\text{for all}\,  e\in \mathcal{E}(K).
\ee
The combination of \eqref{H1Dof:1} and \eqref{H1Dof:3} results in 
\be
\label{H1:2b}
{\rd}\bm{v}|_{e}=0~~\text{for all}\, e\in \mathcal{E}(K). 
\ee
This leads to ${\rd}\bm{v}=\lambda_{1}\lambda_{2}\lambda_{3}r$ for some $r\in P_{k-3}(K)$.
Besides, according to \eqref{H1Dof:4}, 
\be
\label{H1:1}
({\rd}{\bm{v}}, q)_{K}=0~~\text{for all}\, q\in P_{k-3}(K).
\ee
Thus $r=0$ and ${\rd}\bm{v}=0$. This and \eqref{H1:2a}--\eqref{H1:2b} guarantee there exists some $p\in (\lambda_{1}\lambda_{2}\lambda_{3})^{2}P_{k-4}(K)$ such that 
\ben
\bm{v}=\bc\, p.
\een
This and \eqref{H1Dof:5} conclude $\bm{v}=0$.
\end{proof}
{According to \eqref{H1:2a}--\eqref{H1:2b}, the vectorial $H^{1}$ conforming finite element space $V_{h}$ is $H^{1}(\rd)$ conforming.}

\begin{remark}
For $k=3$, $V_{h}$ is a piecewise polynomial of degree $4$, which happens to be $\mathcal{P}_{4}\Lambda^{1}$ presented in \cite[Section 7]{MR2249345}. It is
an ingredient in the BGG approach \cite{MR0578996} for the Arnold-Winther elasticity element \cite{MR2269741}.
\end{remark}


\begin{remark}
It is straight forward that the space $\Sigma_{k,\triangle_{2}}$ is a subset of $H({\rd}{\bd},\Omega;\ms)\cap H({\bd},\Omega;\ms)$. Actually, $\Sigma_{k,\triangle_{2}}$ is able to preserve the Hilbert complex
\begin{center}
\begin{tikzpicture}
\coordinate (a2) at (1,-1.8);
\coordinate (ah2) at (1.5,-1.8);
\coordinate (b2) at (2,-1.8);
\coordinate (c2) at (4.4,-1.8);
\coordinate (ch2) at (5.2,-1.8);
\coordinate (d2) at (6.0,-1.8);
\coordinate (e2) at (11,-1.8);
\coordinate (eh2) at (11.7,-1.8);
\coordinate (f2) at (12.4,-1.8);
\coordinate (g2) at (13.6,-1.8);
\coordinate (h2) at (14,-1.8);

\coordinate[label=left:$RT$] (rt2) at (a2);
 \draw[-latex] (a2)--(b2);
 \coordinate[label=right:${H^{1}({{\rd},\Omega};\mathbb{R}^{2})}$] (h32) at (b2);
 \draw[-latex] (c2)--(d2);
 \coordinate[label=right:${H({{\rd}\bd,\Omega};\ms)\cap H({\bd,\Omega};\ms)}$] (h22) at (d2);
  \draw[-latex] (e2)--(f2);
  \coordinate[label=right:$L^{2}(\Omega)$] (l22) at (f2);  
  \draw[-latex] (g2)--(h2);
   \coordinate[label=right:$0$] (012) at (h2);

  \coordinate[label=above:${\subset}$] (in2) at (ah2);
 \coordinate[label=above:${\rm{sym}\bc}$] (sc2) at (ch2);
  \coordinate[label=above:${{\rd}\bd}$] (dd2) at (eh2);

\end{tikzpicture}
\end{center}
in the discrete case {on a contractible domain $\Omega$}.  The commuting diagram in \cite{Chen2020Finite} can also be constructed here.
\end{remark}

Before establishing the exact complex for the finite elements,  the exact complex for bubble function spaces is constructed below.
Define
\begin{align}
\overset{\circ}{V}_{k+1}(K):&=\{{\bm{v}}\in P_{k+1}(K;\mathbb{R}^{2}):\text{all degrees of freedom \eqref{H1Dof:1}--\eqref{H1Dof:3} vanish} \}.\\
\overset{\circ}{\Sigma}_{k}(K):&=\{{\bm{\sigma}}\in P_{k}(K;\ms):\text{all degrees of freedom \eqref{2dDof:1}--\eqref{2dDof:3} vanish} \}.\\
\overset{\circ}{P}_{k-2}(K):&=P_{k-2}(K)/P_{1}(K).
\end{align}

\begin{lemma}
\label{complex:local}
Given $K\in \mt$, it holds
\ben
{\rd}\bd\overset{\circ}{\Sigma}_{k}(K)=\overset{\circ}{P}_{k-2}(K).
\een
\end{lemma}
\begin{proof}
It is straight forward from \eqref{Green} that
\[{\rd}\bd\overset{\circ}{\Sigma}_{k}(K)\subseteq\overset{\circ}{P}_{k-2}(K).\]
It suffices to prove $\overset{\circ}{P}_{k-2}(K)\subseteq {\rd}\bd\overset{\circ}{\Sigma}_{k}(K)$. Actually, if the inclusion does not hold, then there exists some $q\in \overset{\circ}{P}_{k-2}(K)$, and $q\neq 0$, such that
\ben
({\rd}\bd\bta,q)_{K}=0~~\text{for all}\, \bta\in \overset{\circ}{\Sigma}_{k}(K).
\een
Integration by parts as in \eqref{Green}  leads to
\ben
(\bta,\nabla^{2}q)_{K}=0~~ \text{for all}\,\bta\in \overset{\circ}{\Sigma}_{k}(K).
\een
According to the degrees of freedom \eqref{2dDof:1}--\eqref{2dDof:6}, there exists $\bta\in\overset{\circ}{\Sigma}_{k}(K)$, such that $(\bta,\nabla^{2}q)_{K}\neq0$ as long as $\nabla^{2}q\neq 0$. Hence 
\ben
\nabla^{2}q=0.
\een
This implies $q\in P_{1}(K)$. The contradiction occurs. This concludes the proof.

\end{proof}
\begin{lemma}
\label{poly-complex}
For any triangle $K$, the polynomial complexes 
\begin{center}
\begin{tikzpicture}
\coordinate (a) at (1,0);
\coordinate (ah) at (1.5,0);
\coordinate (b) at (2,0);
\coordinate (c) at (4.3,0);
\coordinate (ch) at (5.15,0);
\coordinate (d) at (6,0);
\coordinate (e) at (7.8,0);
\coordinate (eh) at (8.9,0);
\coordinate (f) at (10,0);
\coordinate (g) at (11.7,0);
\coordinate (h) at (12.5,0);

\coordinate[label=left:$RT$] (rt) at (a);
\draw[-latex] (a)--(b);
\coordinate[label=right:$P_{k+1}(K;\mathbb{R}^{2})$] (vh) at (b);
\draw[-latex] (c)--(d);
\coordinate[label=right:$P_{k}(K;\ms)$] (sh) at (d);
\draw[-latex] (e)--(f);
\coordinate[label=right:$P_{k-2}(K)$] (ph) at (f);  
\draw[-latex] (g)--(h);
\coordinate[label=right:$0$] (0) at (h);

 \coordinate[label=above:${\subset}$] (in) at (ah);
 \coordinate[label=above:${\rm{sym}\bc}$] (sc) at (ch);
 \coordinate[label=above:${{\rd}\bd}$] (dd) at (eh);

\end{tikzpicture}
\end{center}
and 
\begin{center}
\begin{tikzpicture}
\coordinate (a) at (1,0);
\coordinate (ah) at (1.5,0);
\coordinate (b) at (2,0);
\coordinate (c) at (3.8,0);
\coordinate (ch) at (4.9,0);
\coordinate (d) at (6,0);
\coordinate (e) at (7.5,0);
\coordinate (eh) at (8.6,0);
\coordinate (f) at (10,0);
\coordinate (g) at (11.7,0);
\coordinate (h) at (12.5,0);

\coordinate[label=left:$0$] (rt) at (a);
\draw[-latex] (a)--(b);
\coordinate[label=right:$\overset{\circ}{V}_{k+1}(K)$] (vh) at (b);
\draw[-latex] (c)--(d);
\coordinate[label=right:$\overset{\circ}{\Sigma}_{k}(K)$] (sh) at (d);
\draw[-latex] (e)--(f);
\coordinate[label=right:$\overset{\circ}{P}_{k-2}(K)$] (ph) at (f);  
\draw[-latex] (g)--(h);
\coordinate[label=right:$0$] (0) at (h);

 \coordinate[label=above:${\subset}$] (in) at (ah);
 \coordinate[label=above:${\rm{sym}\bc}$] (sc) at (ch);
 \coordinate[label=above:${{\rd}\bd}$] (dd) at (eh);

\end{tikzpicture}
\end{center}
are exact.
\end{lemma}

\begin{proof}
The first polynomial complex follows directly from \cite[Lemma 3.1]{Chen2020Finite}. {The exactness also follows from the existence of homotopy operators; see \cite{MR4066098,2019FiniteHu}}. To obtain the second complex, let $\bsi:={\rm{sym}}{\bc}{\bm{v}}$ for any $ \bm{v}\in \overset{\circ}{V}_{k+1}(K)$, $\bsi\in \overset{\circ}{\Sigma}_{k}(K)$ needs proving.
According to \eqref{H1Dof:1}, $\bsi(a)=0$ for all $a\in \mathcal{V}(K)$. For $\bm{\phi}\in P_{k-2}(e;\mathbb{R}^{2})$, 
\ben
\begin{split}
(\bsi \mathbf{n},\bm{\phi})_{e}&=(\mathbf{n}^{{\rm{T}}}\bsi \mathbf{n}, \bm{\phi}\cdot\mathbf{n})_{e}+(\mathbf{t}^{{\rm{T}}}\bsi \mathbf{n}, \bm{\phi}\cdot\mathbf{t})_{e}\\
&:=\uppercase\expandafter{\romannumeral1}+\uppercase\expandafter{\romannumeral 2}.
\end{split}
\een
The calculations in  \cite[Lemma 2.2]{Chen2020Finite} lead to some identities 
\begin{align}
&\mathbf{n}^{{\rm{T}}}\bsi\mathbf{n}=\mathbf{n}^{{\rm{T}}}\partial_{t}{\bm{v}}, \label{cl:1}\\
&\mathbf{t}^{{\rm{T}}}\bsi\mathbf{n}=\mathbf{t}^{{\rm{T}}}\partial_{t}\bm{v}-\frac{1}{2}\bd\bm{v}, \label{cl:2}\\
&\bd\bsi\cdot\mathbf{n}=\frac{1}{2}\partial_{t}\bd\bm{v}.\label{cl:3}
\end{align}
Combined with \eqref{H1Dof:1}--\eqref{H1Dof:2}, \eqref{cl:1} leads to $\uppercase\expandafter{\romannumeral1}=0$. Combined with \eqref{H1Dof:1}--\eqref{H1Dof:3},  \eqref{cl:2} leads to $\uppercase\expandafter{\romannumeral2}=0$.
The identity \eqref{cl:3} plus \eqref{H1Dof:1} and \eqref{H1Dof:3} result in  
\ben
\bd\bsi\cdot \mathbf{n}=0~~\text{on each}\, e\in \mathcal{E}(K).
\een
The previous arguments lead to $\bsi\in\overset{\circ}{\Sigma}_{k}(K)$, and ${\rm{sym}}{\bc}\, \overset{\circ}{V}_{k+1}(K)\subset \overset{\circ}{\Sigma}_{k}(K)$.

On the other hand, a direct calculating leads to
\ben
\begin{aligned}
&{\rm{dim}}(\overset{\circ}{\Sigma}_{k}(K))=(k-1)^{2}+\frac{(k-2)(k-3)}{2}-4,\\
&{\rm{dim}}({\rd}{\bd}\overset{\circ}{\Sigma}_{k}(K))=\frac{1}{2}k(k-1)-3={\rm{dim}}(\overset{\circ}{P}_{k-2}(K)),\\
&{\rm{dim}}(\overset{\circ}{V}_{k+1}(K))=(k-2)^{2}-1.
\end{aligned}
\een
These result in
\ben
{\rm{dim}}({\rd}{\bd}\overset{\circ}{\Sigma}_{k}(K))={\rm{dim}}(\overset{\circ}{\Sigma}_{k}(K))-{\rm{dim}}(\overset{\circ}{V}_{k+1}(K)).
\een
Together with Lemma \ref{complex:local}, the exactness of the complex follows.
\end{proof}

\begin{lemma}[\cite{2020Complexes}]
\label{divdivcomplex}
The divdiv Hilbert complex 
\begin{center}
\begin{tikzpicture}
\coordinate (a1) at (1,1.8);
\coordinate (ah1) at (1.5,1.8);
\coordinate (b1) at (2,1.8);
\coordinate (c1) at (4,1.8);
\coordinate (ch1) at (5,1.8);
\coordinate (d1) at (6,1.8);
\coordinate (e1) at (7.8,1.8);
\coordinate (eh1) at (8.9,1.8);
\coordinate (f1) at (10,1.8);
\coordinate (g1) at (11.2,1.8);
\coordinate (h1) at (12,1.8);

\coordinate[label=left:$RT$] (rt1) at (a1);
 \draw[-latex] (a1)--(b1);
 \coordinate[label=right:$H^{3}(\Omega;\mathbb{R}^{2})$] (h3) at (b1);
 \draw[-latex] (c1)--(d1);
 \coordinate[label=right:$H^{2}(\Omega;\ms)$] (h2) at (d1);
  \draw[-latex] (e1)--(f1);
  \coordinate[label=right:$L^{2}(\Omega)$] (l2) at (f1);  
  \draw[-latex] (g1)--(h1);
   \coordinate[label=right:$0$] (01) at (h1);

  \coordinate[label=above:${\subset}$] (in1) at (ah1);
 \coordinate[label=above:${\rm{sym}\bc}$] (sc1) at (ch1);
  \coordinate[label=above:${{\rd}\bd}$] (dd1) at (eh1);
\end{tikzpicture}
\end{center}
is exact {on a contractible domain $\Omega$}.
\end{lemma}

Similarly as \cite[Section 3.3]{Chen2020Finite}, the interpolations with commuting properties can be constructed as follows.
Denote the local nodal interpolation operator based on the degrees of freedom \eqref{2dDof:1}--\eqref{2dDof:6} as $\Pi_{K,\triangle_{2}}:H^{2}(K;\ms)\rightarrow P_{k}(K;\ms)$. For any $\bta\in P_{k}(K;\ms)$,  $\Pi_{K,\triangle_{2}}\bta=\bta$ is easy to verify. For the shape regular mesh $\mt$, 
\be
\label{BB:interr}
~~~~~~~\|\bta-\Pi_{K,\triangle_{2}}\bta\|_{0,K}+h_{K}|\bta-\Pi_{K,\triangle_{2}}\bta|_{1,K}+h_{K}^{2}|\bta-\Pi_{K,\triangle_{2}}\bta|_{2,K}\lesssim h_{K}^{s} |\bta|_{s,K}
\ee
holds for $\bta\in H^{s}(K;\ms)$ with $2\leq s\leq k+1$. Integration by parts leads to 
\be
\label{orth}
{\rd}\bd(\Pi_{K,\triangle_{2}}\bta)=\mathcal{Q}_{k-2}^{K}{\rd}{\bd}\bta~~\text{for all}~~ \bta\in H^{2}(K;\ms).
\ee
Here $\mathcal{Q}_{k-2}^{K}:L^{2}(K)\rightarrow P_{k-2}(K)$ is the $L^{2}$ projection operator. It may be later denoted as $\mathcal{Q}_{k-2,\triangle_{d}}^{K}$, $d=2,3$, to distinguish the dimension of $K$.

Denote the local nodal interpolation operator based on the degrees of freedom \eqref{H1Dof:1}--\eqref{H1Dof:5} as $\tilde{I}_{K}:H^{3}(K;\mathbb{R}^{2})\rightarrow P_{k+1}(K;\mathbb{R}^{2})$. For any ${\bm{v}}\in P_{k+1}(K;\mathbb{R}^{2})$,  $\tilde{I}_{K}{\bm{v}}={\bm{v}}$ is easy to verify. For the shape regular mesh $\mt$, 
\be
\label{H1inter:I}
~~~~~~~\|{\bm{v}}-\tilde{I}_{K}{\bm{v}}\|_{0,K}+h_{K}|{\bm{v}}-\tilde{I}_{K}{\bm{v}}|_{1,K}\lesssim h_{K}^{s} |{\bm{v}}|_{s,K}
\ee
holds for ${\bm{v}}\in H^{s}(K;\mathbb{R}^{2})$ with $3\leq s\leq k+2$. The proof of Lemma \ref{poly-complex} shows
\ben
\Pi_{K,\triangle_{2}}({\rm{sym}}\bc\, \bm{v})-{\rm{sym}}\bc(\tilde{I}_{K}{\bm{v}})\in \overset{\circ}{\Sigma}_{k}(K).
\een
Hence, according to Lemma \ref{poly-complex}, there exists $\tilde{\bm{v}}\in\overset{\circ}{V}_{k+1}(K) $ such that
\begin{align}
{\rm{sym}}\bc\tilde{\bm{v}}=\Pi_{K,\triangle_{2}}({\rm{sym}}\bc{\bm{v}})-{\rm{sym}}\bc(\tilde{I}_{K}\bm{v}),\\
\|\tilde{\bm{v}}\|_{0,K}\lesssim h_{K} \|\Pi_{K,\triangle_{2}}({\rm{sym}}\bc{\bm{v}})-{\rm{sym}}\bc(\tilde{I}_{K}\bm{v})\|_{0,K}.
\end{align}
Let $I_{K}\bm{v}:=\tilde{I}_{K}\bm{v}+\tilde{\bm{v}}$. It is also easy to verify $I_{K}\bm{v}=\bm{v}$ for any $\bm{v}\in P_{k+1}(K;\mathbb{R}^{2})$, and
\be
{\rm{sym}}\bc(I_{K}\bm{v})=\Pi_{K,\triangle_{2}}({\rm{sym}}\bc\bm{v})~~\text{for all}~~ \bm{v}\in H^{3}(K;\mathbb{R}^{2}).
\ee
It follows from \eqref{BB:interr} and \eqref{H1inter:I} that
\be
\label{H1inter:II}
~~~~~~~\|{\bm{v}}-{I}_{K}{\bm{v}}\|_{0,K}+h_{K}|{\bm{v}}-{I}_{K}{\bm{v}}|_{1,K}\lesssim h_{K}^{s} |{\bm{v}}|_{s,K}
\ee
with $3\leq s\leq k+2$.

For each $K\in \mt$, let $I_{h}:H^{3}(\Omega;\mathbb{R}^{2})\rightarrow V_{h}$ be defined by $(I_{h}\bm{v})|_{K}:=I_{K}(\bm{v}|_{K})$,  and $\Pi_{h,\triangle_{2}}:H^{2}(\Omega;\ms)\rightarrow \Sigma_{k,\triangle_{2}}$ is defined by $(\Pi_{h,\triangle_{2}}\bta)|_{K}:=\Pi_{K,\triangle_{2}}(\bta|_{K})$, as well as
$\mathcal{Q}_{h,\triangle_{2}}:L^{2}(\Omega)\rightarrow P_{h,\triangle_{2}}$ is defined by $(\mathcal{Q}_{h,\triangle_{2}}q)|_{K}:=\mathcal{Q}_{k-2,\triangle_{2}}^{K}(q|_{K})$.

It follows immediately
\begin{align}
{\rd}{\bd} (\Pi_{h,\triangle_{2}}\bta)=\mathcal{Q}_{h,\triangle_{2}}{\rd}\bd\bta~~\text{for all}~~ \bta\in H^{2}(\Omega;\ms),\label{exchage}\\
{\rm{sym}}\bc(I_{h}\bm{v})=\Pi_{h,\triangle_{2}}({\rm{sym}}\bc\bm{v})~~\text{for all}~~ \bm{v}\in H^{3}(\Omega;\mathbb{R}^{2}).\label{exchage:2}
\end{align}

\begin{lemma}
\label{exactFem}
The {sequence}
\begin{center}
\begin{tikzpicture}
\coordinate (a) at (1,0);
\coordinate (ah) at (1.5,0);
\coordinate (b) at (2,0);
\coordinate (c) at (2.6,0);
\coordinate (ch) at (3.8,0);
\coordinate (d) at (5,0);
\coordinate (e) at (6,0);
\coordinate (eh) at (7,0);
\coordinate (f) at (8,0);
\coordinate (g) at (9,0);
\coordinate (h) at (10,0);

\coordinate[label=left:$RT$] (rt) at (a);
 \draw[-latex] (a)--(b);
 \coordinate[label=right:$V_{h}$] (vh) at (b);
 \draw[-latex] (c)--(d);
 \coordinate[label=right:$\Sigma_{k,\triangle_{2}}$] (sh) at (d);
  \draw[-latex] (e)--(f);
  \coordinate[label=right:$P_{h,\triangle_{2}}$] (ph) at (f);  
  \draw[-latex] (g)--(h);
   \coordinate[label=right:$0$] (0) at (h);

  \coordinate[label=above:${\subset}$] (in) at (ah);
 \coordinate[label=above:${\rm{sym}\bc}$] (sc) at (ch);
  \coordinate[label=above:${{\rd}\bd}$] (dd) at (eh);
\end{tikzpicture}
\end{center}
{is a complex, which is exact on contractible domains}. 
\end{lemma}
\begin{proof}
It is straight forward that 
\ben
{\rd}\bd \Sigma_{k,\triangle_{2}}\subseteq P_{h,\triangle_{2}}.
\een
To obtain $P_{h,\triangle_{2}}={\rd}\bd \Sigma_{k,\triangle_{2}}$, it suffices to prove $P_{h,\triangle_{2}}\subseteq {\rd}\bd \Sigma_{k,\triangle_{2}}$. If the inclusion does not hold, then there exists some $q\in P_{h,\triangle_{2}}$, and $q\neq 0$, such that
\ben
({\rd}\bd \bta, q)_{\Omega}=0~~\text{for all}\, \bta\in  \Sigma_{k,\triangle_{2}}.
\een
Lemma \ref{complex:local} shows, {$q|_{K}\in P_{1}(K)$ for all $K\in \mt$. }
Let the only nonzero degrees of freedom of $\bta$ be
\[(\bd\bta\cdot\mn_{e},[q]_{e})_{e}=([q]_{e},[q]_{e})_{e}~~\text{on some}~~e\in \me_{h}.\]
 Integration by parts leads to 
\ben
0=({\rd}\bd \bta, q)_{\Omega}=\sum_{K\in \mt}({\rd}\bta\cdot\mathbf{n}, q)_{\partial K}=\|[q]_{e}\|_{0, e}^{2}. \\
\een
This shows $[q]_{e}=0$. The arbitrariness of the choice of $e\in \me_{h}$ leads to $q=0$. The contradiction occurs. 

In addition, \eqref{exchage:2} implies
 \ben
 {\rm{sym}}\bc \, V_{h}\subseteq\Sigma_{h,\triangle_{2}}.
 \een

By counting the dimensions, 
\ben
\begin{aligned}
&{\rm{dim}}\, \Sigma_{h,\triangle_{2}}=3\#\mathcal{V}_{h}+(3k-2)\#\mathcal{E}_{h}+\frac{3}{2}k(k-3)\# \mathcal{T}_{h}.\\
&{\rm{dim}}\, {\rm{sym}}\bc\, V_{h}=6\#\mathcal{V}_{h}+(3k-5)\#\mathcal{E}_{h}+(k-1)(k-3)\# \mathcal{T}_{h}-3.\\
&{\rm{dim}}\, {\rd}\bd\, \Sigma_{h,\triangle_{2}}=\frac{1}{2}k(k-1)\#\mathcal{T}_{h}.\\
\end{aligned}
\een
Here $\# \mathcal{S}$ is the number of the elements in the finite set $S$. According to the Euler's formula $\#\mathcal{E}_{h}+1=\#\mathcal{V}_{h}+\#\mathcal{T}_{h}$, 
\begin{align}
{\rm{dim}}\, \Sigma_{h,\triangle_{2}}=&{\rm{dim}}\, {\rm{sym}}\bc\, V_{h}+{\rm{dim}}\, {\rd}\bd\, \Sigma_{h,\triangle_{2}}.
\end{align}
This concludes that the complex is exact.
\end{proof}

\subsection{The construction of the conforming elements on tetrahedral grids}
In this subsection, $\Omega$ is a bounded polyhedron in $\mathbb{R}^{3}$. Given a tetrahedron $K\in \mt$, the finite element shape functions are formed by $P_{k}(K; \ms)$, $k\geq 3$. Some results and notation are introduced here for ensuing use. 
\begin{lemma}[\cite{MfemB}]
\label{3d:6}
Suppose $K\in \mt$, ${\bm{\psi}}\in P_{k}(K;\mathbf{R}^{3})$ satisfies $\bd{\bm{\psi}}=0$, and ${\bm{\psi}}\cdot\mathbf{n}\, |_{\partial K}=0$. Then there exists some $\bm{\vartheta}\in W_{k+1}(K;\mathbb{R}^{3})$ such that 
\ben
\bm{\psi}=\mathbf{curl}\, \bm{\vartheta},
\een
where $W_{k+1}(K;\mathbb{R}^{3})$ is defined by
\ben
W_{k+1}(K;\mathbb{R}^{3}):=\{\bm{\phi}\in P_{k+1}(K; \mathbb{R}^{3}):\, \bm{\phi}\times \mathbf{n}|_{F}=0  ~~\text{for all }~~ F\in \mathcal{F}(K)\}.
\een
\end{lemma}

Define 
\ben
M_{k+2}(K;\ms):=\{\bta\in P_{k+2}(K;\ms):\, (I-\mathbf{n}\mathbf{n}^{{\rm{T}}})\bta (I-\mathbf{n}\mathbf{n}^{{\rm{T}}})|_{F}=\Lambda_{F}(\tau)=0~~\text{for all} ~~F\in \mathcal{F}(K)\}
\een
with $\Lambda_{F}(\bta): F\rightarrow (I-\mathbf{n}\mathbf{n}^{{\rm{T}}}){\ms}(I-\mathbf{n}\mathbf{n}^{{\rm{T}}})$ being defined by 
\ben
\Lambda_{F}(\bta)=(I-\mathbf{n}\mathbf{n}^{{\rm{T}}})(2\varepsilon(\bta\mn)-\partial_{n}\bta) (I-\mathbf{n}\mathbf{n}^{{\rm{T}}}).
\een
Here $F$ is a plane with unit normal $\mn$.

\begin{lemma} [\cite{MR2398766}] 
\label{3d:7}
Suppose $K\in \mt$ is a tetrahedron, $\bta\in P_{k}(K;\ms)$ satisfies $\bd\bta=0$, and $\bta\mathbf{n}|_{\partial K}=0$. Then there exists some $\bm{\zeta} \in  M_{k+2}(K;\ms)$ such that

$$\bta={\mathbf{curl}\mathbf{curl}^{*}}{\bm{\zeta}}.$$ 

\end{lemma}

In addition, define
\be
\label{tem:def:wk}
\mathcal{W}_{k-1}(K; \mathbb{R}^{3}):=\bc\, W_{k}(K;\mathbb{R}^{3})/RM_{\triangle_{3}}(K),
\ee
and
\be
\label{m:def:1}
\mathcal{M}_{k}(K;\ms):=\bc \bc^{*}\, M_{k+2}(K;\ms).
\ee

The degrees of freedom are
\begin{align}
 \bsi (a)&~~ \text{for all}\, a\in \mathcal{V}(K); \label{3dDof:1}\\
(\mathbf{t}_{e}^{{\rm{T}}}\bsi\mathbf{n}_{j}, q)_{e},\, (\mathbf{n}_{i}^{{\rm{T}}}\bsi\mathbf{n}_{j},q)_{e}&~~1\leq i, j \leq 2,  \text{for all}\, q \in P_{{k-2}}(e),\, e\in \mathcal{E}(K); \label{3dDof:2}\\
 (\bsi \mathbf{n}, \bm{\phi})_{F}&~~  \text{for all}\,\bm{\phi} \in P_{{k-3}}(F; \mathbb{R}^{3}), \, F\in \mathcal{F}(K);\label{3dDof:3} \\
 (\bd \bsi\cdot \mathbf{n}, q)_{F}&~~  \text{for all}\, q \in P_{{k-1}}(F), \, F\in \mathcal{F}(K); \label{3dDof:4}\\
 (\bsi, \nabla^{2} q)_{K}&~~  \text{for all}\, q \in P_{{k-2}}(K); \label{3dDof:5}\\
 (\bsi, \nabla  \bm{\phi})_{K}&~~ \text{for all}\, \bm{\phi} \in \mathcal{W}_{k-1}(K; \mathbb{R}^{3}); \label{3dDof:6} \\
  (\bsi,\bm{\tau})_{K}&~~  \text{for all}\, \bm{\tau}\in \mathcal{M}_{k}(K; \ms); \label{3dDof:7}
\end{align}

The degrees of freedom \eqref{3dDof:1}--\eqref{3dDof:3} are exactly the characterization of the continuity of $H(\bd; \ms)$ in \cite{Hu2014A, Hu2015FINITE}, and the continuity of \eqref{3dDof:4} across each interior face leads to $\bd\bsi\in H({\rd}; \mathbb{R}^{3})$.


The global conforming finite element space is defined by
\be
\label{3d:space}
\begin{split}
\Sigma_{k,\triangle_{3}}:=&\{\bm{\tau}\in H({\rd}{\bd},\Omega;\ms):\, \bm{\tau}|_{K}\in P_{k}(K;\ms) ~~ \text{for all}~~ K\in\mt, \\
& \text{all the degrees of freedom \eqref{3dDof:1}--\eqref{3dDof:7} are single-valued}\}.
\end{split}
\ee

\begin{theorem}
The degrees of freedom \eqref{3dDof:1}--\eqref{3dDof:7} uniquely determine a polynomial of $P_{k}(K; \ms)$ defined in \eqref{3d:space}.
\end{theorem}
\begin{proof}
{Note that given any $\bm{\psi}\in W_{k}(K;\mathbb{R}^{3})$, $\bm{\psi}\times\mn|_{\partial K}=0$ and an integration by parts lead to $(\bc\bm\psi, \bm{u})_{K}=(\bm{\psi},\bc \bm{u})_{K}$ for any $\bm{u}\in RM_{\triangle_{3}}(K)$. This, \eqref{tem:def:wk}, and $\bc RM_{\triangle_{3}}(K)=P_{0}(K;\mathbb{R}^{3})$ imply} 
\ben
\begin{split}
{\rm{dim}}\,  \mathcal{W}_{k-1}(K;\mathbb{R}^{3})&={\rm{dim}}\,  \bc\, {W}_{k}(K;\mathbb{R}^{3})-3\\
&={\rm{dim}}\,  {W}_{k}(K;\mathbb{R}^{3})-\text{dim}\,  \nabla \overset{\circ}{P}_{k+1}(K)-3\\
&={\rm{dim}}\,  {RT}_{k-3}(K;\mathbb{R}^{3})-{\rm{dim}}\,  P_{k-3}(K)-3\\
&=\frac{2k^{3}-3k^{2}-5k-12}{6}.\\
\end{split}
\een
The number of degrees of freedom \eqref{3dDof:7} {reads \cite[Theorem 7.2]{MR2398766}}
\ben
{\rm{dim}}\,  \mathcal{M}_{k}(K;\ms)=\frac{k^{3}-3k^{2}-4k+12}{2}.
\een
The remainder degrees of freedom can be counted easily. Thus the number of all the degrees of freedom \eqref{3dDof:1}--\eqref{3dDof:7} is
\ben
(k+1)(k+2)(k+3),
\een
which equals to ${\rm{dim}} P_{k}(K;\ms)$. 

Suppose $\bsi\in P_{k}(K;\ms)$ and all the degrees of freedom \eqref{3dDof:1}--\eqref{3dDof:7} are zero. Then the unisolvence for $P_{k}(K;\ms)$ follows from $\bsi=0$. For $v\in P_{k-2}(K)$, integration by parts and the zero degrees of freedom \eqref{3dDof:3}--\eqref{3dDof:5} lead to 
\ben
\begin{split}
({\rd}\bd \bsi, v)_{K}&=(\bsi, \nabla^{2}v)_{K}-\sum_{F\in\mathcal{F}(K)}(\bsi \mathbf{n}, \nabla v)_{F}+\sum_{F\in\mathcal{F}(K)}({\rd} \bsi\cdot \mathbf{n}, v)_{F}=0.
\end{split}
\een
Therefore, 
\be
\label{3d:t1}
\rd \bd\bsi=0.
\ee
This, Lemma \ref{3d:6}, \eqref{3dDof:4} and \eqref{3d:t1} ensure there exists a function $\bm{\psi}\in W_{k}(K;\mathbb{R}^{3})$ such that
\ben
\bd \bsi=\bc \bm{\psi}.
\een
Furthermore, for all $\bm{\vartheta}\in W_{k}(K;\mathbb{R}^{3})$ with $\mathbf{curl}{\bm{\vartheta}}\bot RM_{\triangle_{3}}(K)$,  \eqref{3dDof:6} and  \eqref{3dDof:3} result in
\be
\label{3d:t2}
\begin{split}
(\bc \bm{\psi}, \bc\bm{\vartheta})_{K}&=(\bd \bsi, \bc\bm{\vartheta})_{K}\\
&=-(\bsi, \nabla \bc\bm{\vartheta})_{K}+\sum_{F\in \mathcal{F}(K)}(\bsi\mathbf{n}, \bc\bm{\vartheta})_{F}=0.
\end{split}
\ee
On the other hand, \eqref{3dDof:1}--\eqref{3dDof:3} lead to the following orthogonality,
\be
\label{3d:orth}
(\bd \bsi, \mathbf{v})_{K}=0~~\text{for all}~~ \mathbf{v}\in RM_{\triangle_{3}}(K).
\ee
This and \eqref{3d:t2} prove ${\bm{\psi}}=\mathbf{0}$. Hence $\bd\bsi=\mathbf{0}$. Furthermore, \eqref{3dDof:1}--\eqref{3dDof:3} lead to $\bsi\mathbf{n}=0$ on $\partial K$. According to Lemma \ref{3d:7},  $\bd\bsi=\mathbf{0}$ entails the relation $\bsi=\bc\bc^{*}\bm{\varphi}$ for some $\bm{\varphi}\in M_{k+2}(K;\ms)$. Consequently, \eqref{3dDof:7} concludes $\bsi=\mathbf{0}$ . 
\end{proof}

\begin{remark}
{Alternatively, $\mathcal{M}_{k}(K;\ms)$ from \eqref{m:def:1} can be defined by 
\ben
\mathcal{M}_{k}(K;\ms):=\{\bsi\in P_{k}(K;\ms):\, \bd\bsi=0, \bsi\mn=0\}.
\een
The number of the basis of the $H(\bd;\ms)-P_{k}$ bubble function space $\Sigma_{K,b}:=\{\bsi\in P_{k}(K;\ms):\bsi\mn=0\}$ introduced in \cite[(2.9)]{MR3301063} is ${\rm{dim}}P_{k-2}(K;\ms)$, and the range of $\bd{\Sigma_{K,b}}$ is $P_{k-1}/RM$. Furthermore, restricted to the bubble functions, the adjoint of $\bd$ operator is $-\bm{\varepsilon}$. The dimension of $\mathcal{M}_{k}(K;\ms)$ can also be derived by the subtraction of the dimension of the range of $\bm{\varepsilon}(P_{k-1}/RM)$ from ${\rm{dim}}P_{k-2}(K;\ms)$, which reads
\be
{\rm{dim}}\mathcal{M}_{k}(K;\ms)=6\, {\rm{dim}}P_{k-2}(K)-3\, {\rm{dim}}P_{k-1}(K)+6=\frac{k^{3}-3k^{2}-4k+12}{2}.
\ee
In addition, the basis functions of space $\mathcal{M}_{k}(K;\ms)$ can be constructed by those bubbles in \cite{MR3301063}.}
\end{remark}

 \begin{remark}
The continuity of $\bd\bsi\cdot\mathbf{n}$ can be replaced by $\partial_{n}(\mathbf{n}^{{\rm{T}}}\bsi\mathbf{n})$. However, different from two dimensions, the replacement can not be done for the interpolation of the degrees of freedom \eqref{3dDof:4}.  In fact, for any $v\in P_{k-1}(F)$, 
 \be
 \label{3d:re}
 \begin{split}
 (\bd\bsi\cdot\mathbf{n},v)_{F}&=(\bd_{F}(\bsi\mathbf{n}),v)_{F}+(\partial_{\mathbf{n}}(\mathbf{n}^{{\rm{T}}}\bsi\mathbf{n}),v)_{F}\\
 &=-(\bsi\mathbf{n},\nabla_{F}v)_{F}+(\mathbf{n}_{\partial F}^{{\rm{T}}}\bsi\mathbf{n},v)_{\partial F}+(\partial_{\mathbf{n}}(\mathbf{n}^{{\rm{T}}}\bsi\mathbf{n}),v)_{F}.\\
 \end{split}
 \ee
 Here $\bd_{F}(\bsi\mathbf{n}):=(\mathbf{n}\times \nabla)\cdot(\mathbf{n}\times(\bsi\mathbf{n}))$, and $\nabla_{F}v:=(\mathbf{n}\times \nabla v)\times\mathbf{n}$.
 The first two terms of \eqref{3d:re} are not any of the degrees of freedom defined in \eqref{3dDof:1}--\eqref{3dDof:7}.
\end{remark}

\section{Mixed finite element methods} 
Recall that the dimension $d$ in this paper is either $2$ or $3$. This section exploits the space $P_{h,\triangle_{d}}$ and the anterior $H({\rd}{\bd},\Omega;\ms)$ conforming finite element spaces $\Sigma_{k,\triangle_{d}}$, $d=2, 3$, to discretize the biharmonic equation. The mixed finite element approximation for \eqref{ConV} is to find $\bsi_{h}\in\Sigma_{k,\triangle_{d}}$, and $u_{h}\in P_{h,\triangle_{d}}$ such that 
\be
\label{DiscreV}
\begin{aligned}
(\bsi_{h},\bta_{h})+({\rd }{\bd}\, \bta_{h}, u_{h}) =& 0~~& \text{for all}~~ \bta_{h}\in \Sigma_{k,\triangle_{d}},\\
({\rd}\, {\bd}\bsi_{h}, v_{h}) = &-(f, v_{h})~~&\text{for all}~~ v_{h}\in P_{h,\triangle_{d}}.
\end{aligned}
\ee

\subsection{BB condition} 
In this subsection, the discrete inf-sup condition is proved to obtain the well-posedness of the mixed finite problem \eqref{DiscreV}.
Define $T(\mathbb{X}):=\{K\in\mt: K\cap \mathbb{X}\neq \emptyset\}$ and $N(T(\mathbb{X})):=\# T(\mathbb{X})$ with $\mathbb{X}$ being a vertice $a\in \mathcal{V}_{h}$ or an edge $e\in \mathcal{E}_{h}$.
The proof of the BB condition is based on a quasi-interpolation $\widetilde{\Pi}_{h,\triangle_{d}}$ with $d=2,3$. 

Recall the $L^{2}$ projection $\mathcal{Q}_{k,\triangle_{d}}^{K}$ onto $P_{k,\triangle_{d}}(K)$.
When $d=2$, define $\widetilde{\Pi}_{h,\triangle_{2}}:H^{1}(\Omega;\ms)\cap\{\bta\in L^{2}(\Omega;\ms):\bd\bta\in H^{1}(\Omega;\mathbb{R}^{2})\}\rightarrow \Sigma_{k,\triangle_{2}}$ as follows: for any $\bta\in H^{1}(\Omega;\ms)\cap\{\bta\in L^{2}(\Omega;\ms):\bd\bta\in H^{1}(\Omega;\mathbb{R}^{2})\}$, 
\ben
\begin{aligned}
&\widetilde{\Pi}_{h,\triangle_{2}}\bta(a)=\frac{1}{N(T(a))}\sum_{K^{\prime}\in T(a)} (\mathcal{Q}_{k,\triangle_{2}}^{K^{\prime}}\bta) (a), \\
& ((\widetilde{\Pi}_{h,\triangle_{2}}\bta) \mathbf{n}, \bm{\phi})_{e}=(\bta \mathbf{n}, \bm{\phi})_{e}~~\text{for all}\, \bm{\phi} \in P_{{k-2}}(e; \mathbb{R}^{2}),\\
&(\bd (\widetilde{\Pi}_{h,\triangle_{2}}\bta)\cdot \mathbf{n},\,  q)_{e}=(\bd \bta\cdot \mathbf{n},\,  q)_{e}~~ \text{for all}\, q \in P_{{k-1}}(e),\\
&(\widetilde{\Pi}_{h,\triangle_{2}}\bta, \nabla^{2} \,q)_{K}=(\bta, \nabla^{2} \,q)_{K}~~ \text{for all}\, q \in P_{{k-2}}(K), \\
& (\widetilde{\Pi}_{h,\triangle_{2}}\bta, \nabla \bc \, q)_{K}=(\bta, \nabla \bc \, q)_{K}~~\text{for all}\, q \in \lambda_{1}\lambda_{2}\lambda_{3} P_{k-3}(K)/P_{0}(K),  \\
& (\widetilde{\Pi}_{h,\triangle_{2}}\bta, \mathcal{J} q)_{K}=(\bta, \mathcal{J} q)_{K}~~\text{for all}\, q\in(\lambda_{1}\lambda_{2}\lambda_{3})^{2}P_{k-4}(K),
\end{aligned}
\een
for each $a\in \mathcal{V}_{h}$, $e\in \mathcal{E}_{h}$ and $K\in \mt$.  

When $d=3$, define $\widetilde{\Pi}_{h,\triangle_{3}}:H^{1}(\Omega;\ms)\cap\{\bta\in L^{2}(\Omega;\ms):\bd\bta\in H^{1}(\Omega;\mathbb{R}^{3})\}\rightarrow \Sigma_{k,\triangle_{3}}$ as follows: for any $\bta\in H^{1}(\Omega;\ms)\cap\{\bta\in L^{2}(\Omega;\ms):\bd\bta\in H^{1}(\Omega;\mathbb{R}^{3})\}$, 
\ben
\begin{aligned}
&\widetilde{\Pi}_{h,\triangle_{3}}\bta(a)=\frac{1}{N(T(a))}\sum_{K^{\prime}\in T(a)} (\mathcal{Q}_{k,\triangle_{3}}^{K^{\prime}}\bta) (a), \\
&(\mathbf{t}_{e}^{{\rm{T}}}\widetilde{\Pi}_{h,\triangle_{3}}\bta\mathbf{n}_{j}, q)_{e}=\frac{1}{N(T(e))}\sum_{K^{\prime}\in T(e)}(\mathbf{t}_{e}^{{\rm{T}}}(\mathcal{Q}_{k,\triangle_{3}}^{K^{\prime}}\bta)\mathbf{n}_{j}, q)_{e},~~1\leq  j \leq 2,  \text{for all}\, q \in P_{{k-2}}(e),\\
& (\mathbf{n}_{i}^{{\rm{T}}}\widetilde{\Pi}_{h,\triangle_{3}}\bta\mathbf{n}_{j},q)_{e}=\frac{1}{N(T(e))}\sum_{K^{\prime}\in T(e)}(\mathbf{n}_{i}^{{\rm{T}}}(\mathcal{Q}_{k,\triangle_{3}}^{K^{\prime}}\bta)\mathbf{n}_{j},q)_{e}, ~~1\leq i, j \leq 2,  \text{for all}\, q \in P_{{k-2}}(e),\\
&(\widetilde{\Pi}_{h,\triangle_{3}}\bta \mathbf{n}, \bm{\phi})_{F}= (\bta\mathbf{n}, \bm{\phi})_{F}, ~~~~\text{for all}\,\bm{\phi} \in P_{{k-3}}(F; \mathbb{R}^{3}),\\
&(\bd \widetilde{\Pi}_{h,\triangle_{3}}\bta\cdot \mathbf{n}, q)_{F}=(\bd \bta\cdot \mathbf{n}, q)_{F}, ~~~~ \text{for all}\, q \in P_{{k-1}}(F),\\
&(\widetilde{\Pi}_{h,\triangle_{3}}\bta, \nabla^{2} q)_{K}= (\bta, \nabla^{2} q)_{K},~~~~  \text{for all}\, q \in P_{{k-2}}(K),\\
& (\widetilde{\Pi}_{h,\triangle_{3}}\bta, \nabla  \bm{\phi})_{K}=(\bta, \nabla  \bm{\phi})_{K}, ~~~~\text{for all}\, \bm{\phi} \in \mathcal{W}_{k-1}(K; \mathbb{R}^{3}), \\
& (\widetilde{\Pi}_{h,\triangle_{3}}\bta,{\bm{\psi}})_{K}=(\bta,{\bm{\psi}})_{K}, ~~~~\text{for all}\, {\bm{\psi}}\in \mathcal{M}_{k}(K; \ms),
\end{aligned}
\een
for each $a\in \mathcal{V}_{h}$, $e\in \mathcal{E}_{h}$, $F\in \mathcal{F}_{h}$, as well as $K\in \mt$.

\begin{theorem}
\label{2d:BB}
Assume the triangulation $\mt$ is shape regular. There exists a constant $\beta$ independent of $h$ such that the following BB  condition holds, 
\be
\label{2d:bb1}
\inf_{v_{h}\in P_{h,\triangle_{d}}}\sup_{\bta_{h}\in\Sigma_{k,\triangle_{d}}}\frac{({\rd}{\mathbf{div}}\bta_{h}, v_{h})}{\|\bta_{h}\|_{H({\rd}{\bd})}\|v_{h}\|_{0}}\geq \beta>0.
\ee
Furthermore, the stability for \eqref{DiscreV} is obtained.
\end{theorem}
\begin{proof}
For any $v_{h}\in P_{h,\triangle_{d}}$, according to \cite{Raviart1986Finite}, there exists some ${\bm{\phi}}\in H^{1}(\Omega;\mathbb{R}^{d})$, such that ${\rd}{\bm{\phi}}=v_{h}$, and $\|{\bm{\phi}}\|_{1}\lesssim \|v_{h}\|_{0}$. There exists some $\bta_{0}\in H^{1}(\Omega;\ms)$, such that ${\bd\bta_{0}}={\bm{\phi}}$, and $\|{\bta_{0}}\|_{1}\lesssim \|{\bm{\phi}}\|_{0}$. 

For any $q\in P_{h,\triangle_{d}}$, integration by parts leads to 
\ben
\begin{split}
({\rd}{\bd}\widetilde{\Pi}_{h,\triangle_{d}}\bta, q)&=(\widetilde{\Pi}_{h,\triangle_{d}}\bta,\nabla^{2}q)-\sum_{K\in \mt}(\widetilde{\Pi}_{h,\triangle_{d}}\bta\cdot \mathbf{n},\nabla q)_{\partial K}+\sum_{K\in\mt}(\bd(\widetilde{\Pi}_{h,\triangle_{d}}\bta)\cdot\mathbf{n}, q)_{\partial{K}}\\
&=(\bta,\nabla^{2}q)-\sum_{K\in \mt}(\bta\cdot \mathbf{n},\nabla q)_{\partial K}+\sum_{K\in \mt}(\bd\bta\cdot\mathbf{n},q)_{\partial K}\\
&=({\rd}{\bd}\bta, q).
\end{split}
\een
This implies
\be
\label{BB:exchange}
{\rd}{\bd}\widetilde{\Pi}_{h,\triangle_{d}}\bta=\mathcal{Q}_{h,\triangle_{d}}({\rd}{\bd}\bta).
\ee
The estimates 
\be
\label{err:tilde}
\|\bta-\Pi_{h,\triangle_{d}}\bta\|_{i}\lesssim h^{s-i}|\bta|_{s}+h^{s+1-i}\|\bd\bta\|_{s}, ~~s\geq 1,~~i=0,1
\ee
follow by standard techniques.

Due to \eqref{BB:exchange}--\eqref{err:tilde}, it holds
\ben
\|{\rd}\bd\widetilde{\Pi}_{h,\triangle_{d}}\bta_{0}\|_{0}=\|\mathcal{Q}_{h,\triangle_{d}}{\rd}\bd\bta_{0}\|_{0}=\|\mathcal{Q}_{h,\triangle_{d}}{\rd}{\bm{\phi}}\|_{0}\lesssim \|{\bm{\phi}}\|_{1}\lesssim \|v_{h}\|_{0}.
\een
Thus $ \|\widetilde{\Pi}_{h,\triangle_{d}}\bta_{0}\|_{H({\rd}{\bd})}\lesssim \|v_{h}\|_{0}$. The replacement $\bta_{h}=\widetilde{\Pi}_{h,\triangle_{d}}\bta_{0}$ proves the BB condition \eqref{2d:bb1}.

Additionally, by the Babu$\check{\rm{s}}$ka Brezzi theory \cite{brezzi2012mixed,MfemB}, for any $\widetilde{\bta}_{h}\in \Sigma_{k,\triangle_{d}}$ and $\widetilde{v}_{h}\in P_{h,\triangle_{d}}$,
\be
\label{2d:error}
\begin{split}
&\|\widetilde{\bta}_{h}\|_{H({\rd}{\bd})}+\|\widetilde{v}_{h}\|_{0}\\
&\lesssim \sup_{\substack{\bta_{h}\in \Sigma_{k,\triangle_{d}},\\ v_{h}\in P_{h,\triangle_{d}}}}\frac{(\widetilde{\bta}_{h}, \bta_{h})+({\rd}{\bd}\bta_{h},\widetilde{v}_{h})+({\rd}\bd \widetilde{\bta}_{h}, v_{h})}{\|\bta_{h}\|_{\mathbf{H}({\rd}\bd)}+\|v_{h}\|_{0}}.
\end{split}
\ee
This ensures that the problem \eqref{DiscreV} is well-posed.
\end{proof}
\begin{remark}
The $H^{2}(\Omega;\ms)$ regularity for $\bta$ is required if one employs the interpolation operator $\Pi_{h,\triangle_{2}}$ in \eqref{2d:bb1}. The proof of Theorem \ref{2d:BB} somehow reduces the regularity requirement of the interpolation. Nevertheless, the exactness of the divdiv Hilbert complex in Lemma \ref{divdivcomplex} ensures the existence of  $\bta\in H^{2}(\Omega;\ms)$ for any $v_{h}\in P_{h,\triangle_{2}}$.
\end{remark}


\subsection{Error analysis}
The stability of \eqref{DiscreV} allows the following error estimates.
\begin{theorem}
\label{err:err}
Let $(\bsi,u)\in H({\rd}\mathbf{div},\Omega;\ms)\times L^{2}(\Omega)$ be the solution of \eqref{ConV} and $(\bsi_{h},u_{h})\in\Sigma_{k, \triangle_{d}}\times P_{h,\triangle_{d}} $ be the solution of \eqref{DiscreV}. Assume $\bsi\in H^{k+1}(\Omega;\ms)$, $u\in H^{k-1}(\Omega)$, and $f\in H^{k-1}(\Omega)$, $k\geq 3$. Then

\begin{align}
&\|\bsi-\bsi_{h}\|_{0}\lesssim h^{k+1}|\bsi|_{k+1},\label{err:1}\\
&\|u-u_{h}\|_{0}\lesssim h^{k+1}|\bsi|_{k+1}+h^{k-1}|u|_{k-1},\label{err:2}\\
&\|\bsi-\bsi_{h}\|_{H({\rd}\bd)}\lesssim h^{k+1}|\bsi|_{k+1}+h^{k-1}|f|_{k-1}.\label{err:3}
\end{align}

\end{theorem}
\begin{proof}
Theorem \ref{2d:BB} leads to 
\ben
\begin{split}
&~~~~\|\Pi_{h,\triangle_{d}}\bsi-\bsi_{h}\|_{H({\rd}{\bd})}+\|\mathcal{Q}_{h,\triangle_{d}}u-u_{h}\|_{0} \\
&\lesssim\sup_{\substack{\bta_{h}\in \Sigma_{k,\triangle_{d}},\\ v_{h}\in P_{h,\triangle_{d}}}}\frac{(\Pi_{h,\triangle_{d}}\bsi-\bsi_{h}, \bta_{h})+({\rd}{\bd}\bta_{h},\mathcal{Q}_{h,\triangle_{d}}u-u_{h})+({\rd}\bd(\Pi_{h,\triangle_{d}}\bsi-\bsi_{h}), v_{h})}{\|\bta_{h}\|_{H({\rd}\bd)}+\|v_{h}\|_{0}}.\\
\end{split}
\een
According to \eqref{ConV} and \eqref{DiscreV}, 
\ben
\begin{split}
&(\Pi_{h,\triangle_{d}}\bsi-\bsi_{h}, \bta_{h})+({\rd}{\bd}\bta_{h},\mathcal{Q}_{h,\triangle_{d}}u-u_{h})+({\rd}\bd(\Pi_{h,\triangle_{d}}\bsi-\bsi_{h}), v_{h})\\
&=(\Pi_{h,\triangle_{d}}\bsi-\bsi, \bta_{h}).\\
\end{split}
\een
This shows
\be
\label{err:0}
\|\Pi_{h,\triangle_{d}}\bsi-\bsi_{h}\|_{H({\rd}{\bd})}+\|\mathcal{Q}_{h,\triangle_{d}}u-u_{h}\|_{0}\lesssim \|\Pi_{h,\triangle_{d}}\bsi-\bsi\|_{0}.
\ee
Additionally, together with \eqref{BB:interr} the standard interpolation error estimates for $d=2$, as well as the same results hold for $d=3$, the convergence results  \eqref{err:1}--\eqref{err:3} follow from \eqref{err:0}.  
\end{proof}

\subsection{Superconvergence}
Introduce the space
\ben
H^{2}(\mt):=\{v\in L^{2}(\Omega):v|_{K}\in H^{2}(K)~~ \text{for all}~~ K\in \mt\}.
\een
Define the corresponding mesh-dependent norm in two dimensions,
\ben
|v|_{2,h,\triangle_{2}}^{2}:=\sum_{K\in \mt}|v|_{2,K}^{2}+\sum_{e\in \mathcal{E}_{h}}(h_{e}^{-3}\|[v]\|_{0,e}^{2}+h_{e}^{-1}\|[\nabla v]_{e}\|_{0,e}^{2}),
\een
as well as in three dimensions,
\ben
|v|_{2,h,\triangle_{3}}^{2}:=\sum_{K\in \mt}|v|_{2,K}^{2}+\sum_{F\in \mathcal{F}_{h}}(h_{F}^{-3}\|[v]\|_{0,F}^{2}+h_{F}^{-1}\|[\nabla v]_{F}\|_{0,F}^{2}).
\een

\begin{lemma}
For $d$ being either $2$ or $3$, there exists some constant $\beta>0$ such that the following BB condition regarding to the mesh-dependent norm holds, 
\be
\label{BB:mesh1}
 \sup_{\bta_{h}\in \Sigma_{k,\triangle_{d}}}\frac{({\rd}\bd\bta_{h}, v_{h})}{\|\bta_{h}\|_{0}}\geq \beta |v_{h}|_{2,h,\triangle_{d}}~~\text{for all}~~ v_{h}\in P_{h,\triangle_{d}}.
\ee
\end{lemma}
\begin{proof}
Let $v_{h}\in P_{h,\triangle_{d}}$. For $d=2$, let the degrees of freedom of $\bta_{h}\in \Sigma_{k,\triangle_{2}}$ for each $K\in \mt$ being 
\ben
\begin{aligned}
\bta_{h}(a)&=0~~ &\text{for all}\, a\in \mathcal{V}(K),\\
 (\bta_{h} \mathbf{n}, \bm{\phi})_{e}&=(h_{e}^{-1}[\nabla v_{h}], \bm{\phi})_{e}~~ &\text{for all}\, \bm{\phi} \in P_{{k-2}}(e; \mathbb{R}^{2}), \, e\in \mathcal{E}(K),\\
 (\bd \bta_{h}\cdot \mathbf{n},\,  q)_{e}&=-(h_{e}^{-3}[v_{h}], \bm{\phi})_{e}~~ &\text{for all}\, q \in P_{{k-1}}(e), \, e\in \mathcal{E}(K),\\
 (\bta_{h}, \nabla^{2} \,q)_{K}&=(\nabla^{2} v_{h}, \nabla^{2} \,q)_{K}~~ &\text{for all}\, q \in P_{{k-2}}(K),\\
 (\bta_{h}, \nabla \bc \, q)_{K}&=0~~&\text{for all}\,q \in \lambda_{1}\lambda_{2}\lambda_{3} P_{k-3}(K)/P_{0}(K),\\
  (\bta_{h}, \mathcal{J} q)_{K}&=0~~&\text{for all}\, q\in \overset{\circ}{\mathcal{B}}_{Arg, k+2}(K).
  \end{aligned}
\een
Consider
\ben
\begin{split}
({\rd}{\bd}\bta_{h}, v_{h})&=\sum_{K\in\mt}({\rd}{\bd}\bta_{h}, v_{h})_{K}\\
&=\sum_{K\in\mt}(\bta_{h}, \nabla^{2}v_{h})_{K}-\sum_{e\in \mathcal{E}_{h}}(\bd\bta_{h}\cdot\mathbf{n}, [v_{h}])_{e}+\sum_{e\in \mathcal{E}_{h}}(\bta_{h}\mathbf{n}, [\nabla v_{h}])_{e}\\
&=\sum_{K\in\mt}\|\nabla^{2}v_{h}\|_{0}^{2}+\sum_{e\in \mathcal{E}_{h}}(h_{e}^{-3}\|[v_{h}]\|_{0, e}^{2}+h_{e}^{-1}\|[\nabla v_{h}]\|^{2}_{0, e})\\
&=|v_{h}|_{2,h,\triangle_{2}}^{2}.
\end{split}
\een
The scaling argument leads to 
\ben
\|\bta_{h}\|_{0}\lesssim |v_{h}|_{2,h,\triangle_{2}}.
\een
Therefore, 
\ben
\frac{({\rd}\bd\bta_{h}, v_{h})}{\|\bta_{h}\|_{0}}\gtrsim |v_{h}|_{2,h,\triangle_{2}}.
\een
This proves \eqref{BB:mesh1} in two dimensions.

 When it comes to $d=3$, the same techniques are applied. Let $\bta_{h}\in \Sigma_{k,\triangle_{3}}$ on each $K\in \mt$ with
 \ben
\begin{aligned}
 \bta_{h} (a)&=0~~&\text{for all}\, a\in \mathcal{V}(K),\\
(\mathbf{t}_{e}^{{\rm{T}}}\bsi\mathbf{n}_{j}, q)_{e}&=0, (\mathbf{n}_{i}^{{\rm{T}}}\bsi\mathbf{n}_{j},q)_{e}=0~~&1\leq i, j \leq 2, \text{for all}\, q \in P_{{k-2}}(e),\, e\in \mathcal{E}(K), \\
 (\bta_{h} \mathbf{n}, \bm{\phi})_{F}&=(h_{F}^{-1}[\nabla v_{h}],\bm{\phi} )_{F}~~ &\text{for all}\, \bm{\phi} \in P_{{k-3}}(F; \mathbb{R}^{3}), \, F\in \mathcal{F}(K), \\
 (\bd \bta_{h}\cdot \mathbf{n}, q)_{F}&=-(h_{F}^{-3}[v_{h}], q)_{F}~~ &\text{for all}\, q \in P_{{k-1}}(F), \, F\in \mathcal{F}(K),\\
 (\bta_{h}, \nabla^{2} q)_{K}&=0~~ &\text{for all}\, q \in P_{{k-2}}(K),\\
 (\bta_{h}, \nabla  \bm{\phi})_{K}&=0~~&\text{for all}\, \bm{\phi} \in \mathcal{W}_{k-1}(K; \mathbb{R}^{3}), \\
  (\bta_{h},\bm{\tau})_{K}&=0~~&\text{for all}\, \bm{\tau}\in \mathcal{M}_{k}(K; \ms).
\end{aligned}
\een
This leads to
\ben
({\rd}\bd\bta_{h}, v_{h})=|v_{h}|_{2,h,\triangle_{3}}^{2}.
\een
The scaling argument in this scenario results in 
\ben
\|\bta_{h}\|_{0}\lesssim |v_{h}|_{2,h,\triangle_{3}}.
\een
This proves \eqref{BB:mesh1} in three dimensions.

\end{proof}

Babu$\check{\rm{s}}$ka Brezzi theory \cite{brezzi2012mixed,MfemB} and the BB condition \eqref{BB:mesh1} lead to the following stability results. 
For any $\widetilde{\bta}_{h}\in \Sigma_{k,\triangle_{d}}$ and $\widetilde{v}_{h}\in P_{h,\triangle_{d}}$,
\be
\label{mesh:error}
\begin{split}
&\|\widetilde{\bta}_{h}\|_{0}+|\widetilde{v}_{h}|_{2,h,\triangle_{d}}\\
&\lesssim \sup_{\substack{\bta_{h}\in \Sigma_{k,\triangle_{d}},\\ v_{h}\in P_{h,\triangle_{d}}}}\frac{(\widetilde{\bta}_{h}, \bta_{h})+({\rd}{\bd}\bta_{h},\widetilde{v}_{h})+({\rd}\bd \widetilde{\bta}_{h}, v_{h})}{\|\bta_{h}\|_{0}+|v_{h}|_{2,h,\triangle_{d}}}.
\end{split}
\ee

The stability result \eqref{mesh:error} gives rise to the following superconvergence results.
\begin{theorem}
Suppose $(\bsi_{h},u_{h})\in \Sigma_{k, \triangle_{d}}\times P_{h,\triangle_{d}} $ is the solution of the mixed finite element method \eqref{DiscreV}. Assume $\bsi\in H^{k+1}(\Omega;\ms)$. Then 
\ben
|\mathcal{Q}_{h}u-u_{h}|_{2, h,\triangle_{d}}\lesssim h^{k+1}|\bsi|_{k+1}.
\een
\end{theorem}

\subsection{Postprocessing}
The superconvergence of $|Q_{h}u-u_{h}|_{2, h,\triangle_{d}}$ is used to get a high order approximation of displacement in this subsection. Define $u_{h}^{*}\in P_{k+2}(\mt)$ as follows: for each $K\in \mt$,
\begin{align}
(\nabla^{2}u_{h}^{*},\nabla^{2}q)_{K}&=-(\bsi_{h},\nabla^{2}q)_{K}~~&\text{for all}~~ q\in P_{k+2}(\mt).\\
(u_{h}^{*}, q)_{K}&=(u_{h}, q)_{K}~~&\text{for all}~~ q\in P_{1}(\mt).
\end{align} 

\begin{theorem}
\label{Post}
Suppose $(\bsi_{h},u_{h})\in \Sigma_{h, \triangle_{d}}\times  P_{h,\triangle_{d}}$ is the solution of the mixed finite element method \eqref{DiscreV}. Assume $u\in H^{k+3}(\Omega;\ms)$. Then 
\ben
|u-u_{h}^{*}|_{2, h,\triangle_{d}}\lesssim h^{k+1}|u|_{k+3}.
\een
\end{theorem}
\begin{proof}
The proof of Theorem \ref{Post} is similar as \cite[Theorem 4.4]{Chen2020Finite}, and the details are omitted here.
\end{proof}

\section{Numerical results}
Some numerical results are presented in this section to verify the error analysis and convergence results in previous sections.  
\begin{figure}
\centering
\begin{minipage}[t]{0.5\textwidth}
            \centering          
            \includegraphics[width=0.7\textwidth]{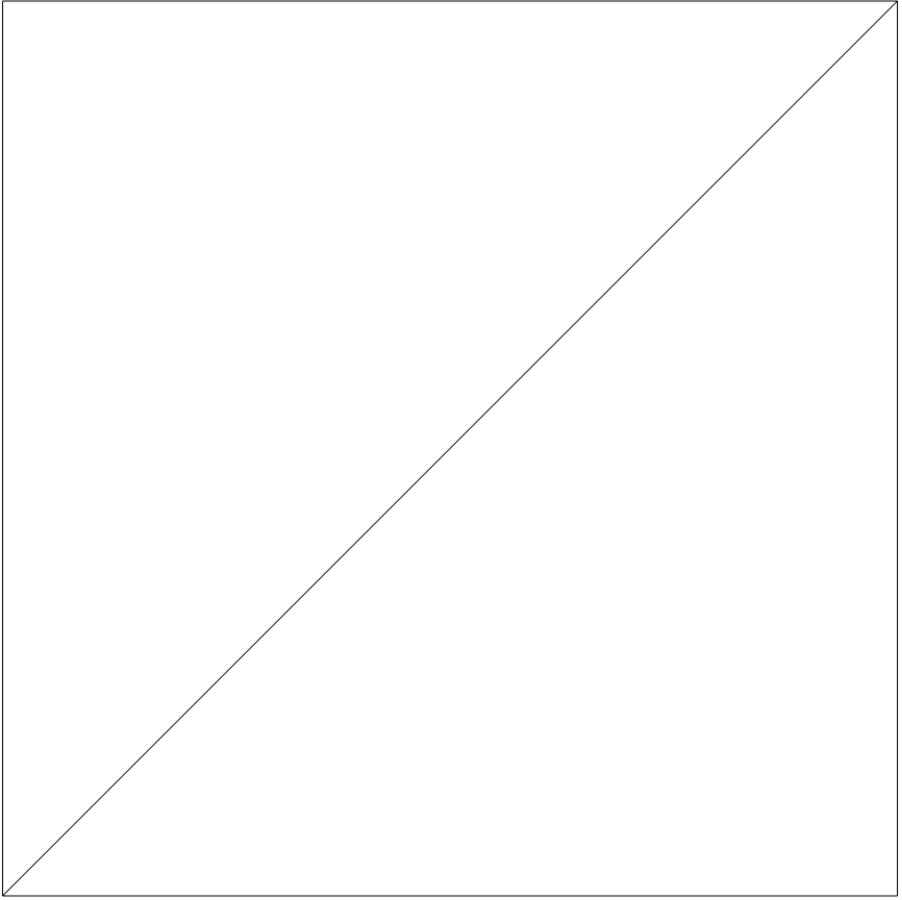}   
            \caption{Initial mesh of the uniform triangulation}
            \label{fig:mesh:a}
           \end{minipage}\begin{minipage}[t]{0.5\textwidth}
            \centering       
            \includegraphics[width=0.7\textwidth]{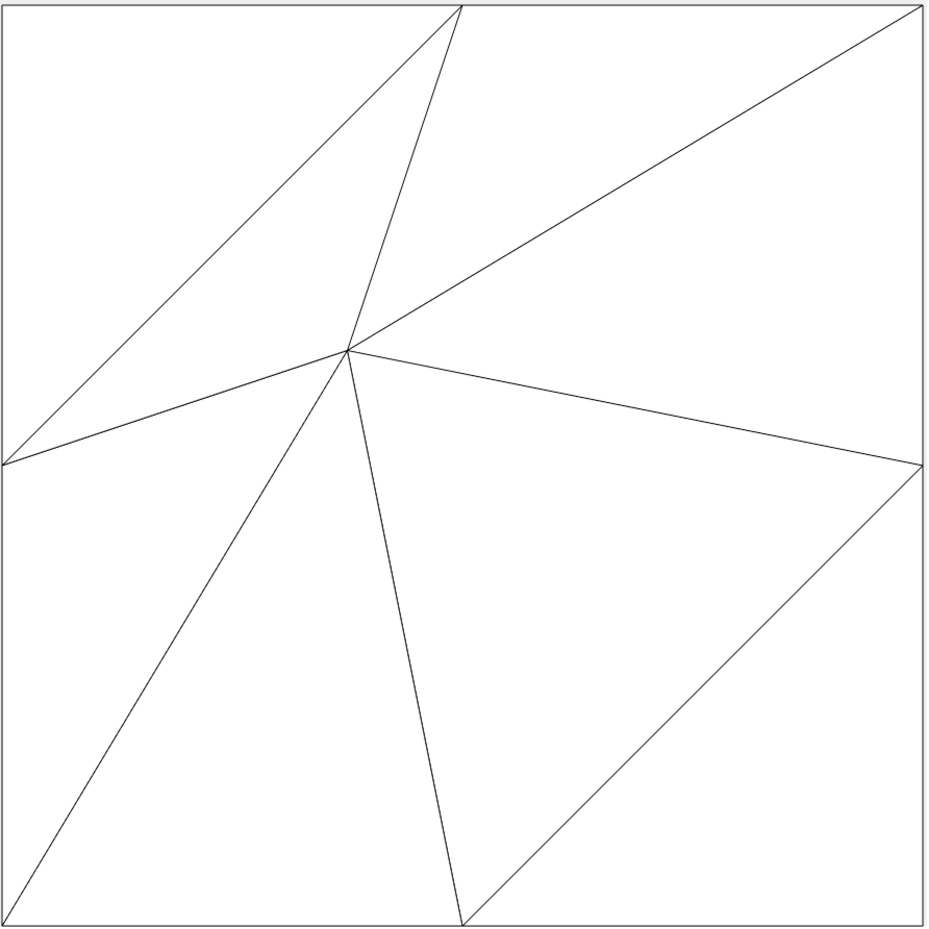}   
            \caption{Initial mesh of the non-uniform triangulation}
            \label{fig:mesh:b}
           \end{minipage}
\begin{minipage}[t]{0.5\textwidth}
            \centering       
            \includegraphics[width=0.7\textwidth]{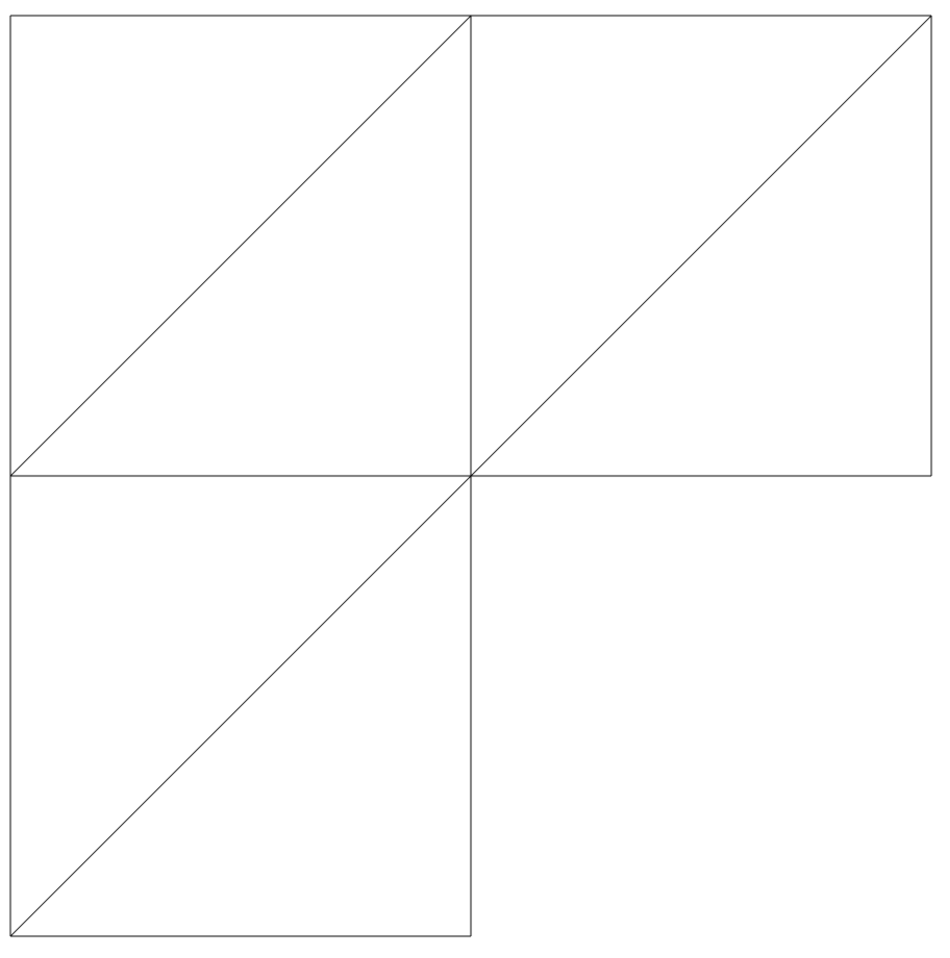}   
           \caption{Initial mesh of the triangulation for the L-shape domain}
            \label{fig:mesh:c}
           \end{minipage}
\end{figure}

\subsection{Example 1}
The computational domain is $\Omega=(0,1)\times (0,1)$ with the homogeneous boundary condition. Load function $f=\Delta^{2} u$ in \eqref{biha} is derived by the exact solution  
$$u(x,y)= x^{2}y^{2}(y - 1)^{2}( 1-x)^{2}.$$ 

Use the $H({\rd}\bd;\ms)$ conforming finite element $\Sigma_{3,\triangle_{2}}$ for $\bsi_{h}$ in problem \eqref{DiscreV}, and the piecewise linear space $P_{h,\triangle_{2}}$ for $u_{h}$. $\mt$ is uniform in this example. The initial mesh is shown in Figure \ref{fig:mesh}. The errors are reported in Table \ref{Ta1}. As shown in Theorem \ref{err:err}, the optimal order of convergence for both $\bsi$ and $u$ is achieved in the computation. Besides, the superconvergence can be observed. The errors $\|Q_{h}u-u_{h}\|_{0}$ and $|Q_{h}u-u_{h}|_{2,h,\triangle_{2}}$ are {fourth} order of convergence, and $|Q_{h}u-u_{h}|_{2,h,\triangle_{2}}$ are {fourth} order higher than the optimal one. In addition, {fourth} order of convergence is achieved for $|u-u_{h}^{*}|_{2,h,\triangle_{2}}$ with the postprocessing solution $u_{h}^{*}$.
\begin{table}[h!]\footnotesize\tabcolsep 16pt
\begin{center}
\caption{The error and the order of convergence on uniform meshes }\label{Ta1}\vspace{-2mm}
\end{center}
\begin{center}
\begin{tabular}{c| c c c c c c }
\toprule  
      &$\|\sigma-\sigma_{h}\|_{0}$    &$h^{n}$     & $\|{\rd}{\bd}(\sigma-\sigma_{h})\|_{0}$   &$h^{n}$  &$\|u-u_{h}\|_{0}$  &$h^{n}$\\
\midrule  
 1 &1.3900e-02     &$-$     &1.5552e+00  &$-$ &9.2528e-04 &$-$ \\
 2  &5.1722e-03  &1.43     &1.0180e+00 & 0.61 &2.8527e-04 &1.70  \\
 3  &4.2279e-04   &3.61    &3.0838e-01   &1.72  &1.0686e-04 &1.42 \\
 4  &2.9243e-05   &3.85    &8.0510e-02   &1.94 &3.0080e-05  &1.83 \\
 5  &1.9079e-06     &3.94   &2.0341e-02  &1.98 &7.7431e-06  &1.96 \\
\toprule
 &$\|Q_{h}u-u_{h}\|_{0}$    &$h^{n}$        &$|Q_{h}u-u_{h}|_{2,h,\triangle_{2}}$ &$h^{n}$  &$|u-u_{h}^{*}|_{2,h,\triangle_{2}}$  &$h^{n}$\\
\midrule
1&3.2820e-04&$-$  &2.8601e-03 &$-$  &2.4999e-02&$-$\\
2 &9.7954e-05  &1.74 &1.8199e-03   &0.65 &5.9936e-03 &2.06\\
3&7.3858e-06   &3.73 &1.7841e-04  &3.35 &5.0637e-04 &3.57\\
4&4.7511e-07   &3.96 &1.3361e-05   &3.74  &3.4275e-05 &3.88\\
 5 &2.9848e-08   &3.99 &9.1053e-07   &3.88 &2.1900e-06 &3.97\\
\bottomrule 
\end{tabular}
\end{center}
\end{table}

\subsection{Example 2}
Compute Example 1 on non-uniform triangulations. The initial mesh is shown in Figure \ref{fig:mesh}. The errors and convergence rates are displayed in Table \ref{Ta:ire}.  The computation shows that the nonuniformity of the mesh does not downgrade approximability. 
\begin{table}[h!]\footnotesize\tabcolsep 16pt
\begin{center}               
\caption{The error and the order of convergence on non-uniform meshes }\label{Ta:ire}\vspace{-2mm}
\end{center}
\begin{center}
\begin{tabular}{c| c c c c c c } 
\toprule
      &$\|\sigma-\sigma_{h}\|_{0}$    &$h^{n}$     & $\|{\rd}{\bd}(\sigma-\sigma_{h})\|_{0}$   &$h^{n}$  &$\|u-u_{h}\|_{0}$  &$h^{n}$\\
\midrule 
 1 &5.7827e-03   &$-$ &1.0671e+00  &$-$ &3.4324e-04  &$-$ \\
 2 &4.4357e-04 &3.70 &3.1255e-01  &1.77  &1.1681e-04    &1.56\\
 3 &3.1802e-05  &3.80 &8.1188e-02   &1.94 &3.3490e-05    &1.80\\
 4 &2.1183e-06   &3.91 &2.0490e-02  &1.99  &8.6457e-06   &1.95 \\
 5 &1.3611e-07   &3.96&5.1346e-03   &2.00 &2.1787e-06   &1.99 \\
\toprule
 &$\|Q_{h}u-u_{h}\|_{0}$    &$h^{n}$        &$|Q_{h}u-u_{h}|_{2,h,\triangle_{2}}$ &$h^{n}$  &$|u-u_{h}^{*}|_{2,h,\triangle_{2}}$  &$h^{n}$\\
\midrule
 1&1.1435e-04 &$-$  &1.9830e-03  &$-$ &6.7217e-03&$-$\\
2 &7.7707e-06 & 3.88 &1.8557e-04  & 3.42&5.3202e-04 &3.66\\
  3&5.1056e-07&3.93 &1.4758e-05 & 3.65 &3.8076e-05 &3.80\\
  4&3.2260e-08  &3.98 &1.0220e-06  &3.85 &2.5243e-06 &3.91\\
 5 &2.0211e-09   &4.00&6.6842e-08   &3.93&1.6174e-07 &3.96\\
\bottomrule
\end{tabular}
\end{center}
\end{table}

\subsection{Example 3}
The L-shape domain $\Omega=(-1,1)\times (-1,1) \backslash ([0,1]\times[-1,0])$. Figure \ref{fig:mesh} shows its initial mesh. Let $\omega:=3\pi/2$, and $\alpha=0.544483736782464$ is a non-characteristic root of $\sin^{2}(\alpha\omega)=\alpha^{2}\sin^{2}(\omega)$ with 
\ben
\begin{split}
g_{\alpha,\omega}(\theta)&=g_{1}(\cos((\alpha-1)\theta)-\cos((\alpha+1)\theta))\\
&-g_{2}(\frac{1}{\alpha-1}\sin((\alpha-1)\theta)-\frac{1}{\alpha+1}\sin((\alpha+1)\theta)),\\
\end{split}
\een
and 
\ben
\begin{aligned}
g_{1}&=\frac{1}{\alpha-1}\sin((\alpha-1)\omega)-\frac{1}{\alpha+1}\sin((\alpha+1)\omega),\\
g_{2}&=\cos((\alpha-1)\omega)-\cos((\alpha+1)\omega).
\end{aligned}
\een
Load function $f=\Delta^{2} u$ in \eqref{biha} is derived by the exact solution
\ben
u(x,y)=(1-x^{2})^{2}(1-y^{2})^{2}(\sqrt{x^{2}+y^{2}})^{1+\alpha}g_{\alpha,\omega}(\theta).
\een

Use the $H({\rd}\bd;\ms)$ conforming finite element $\Sigma_{3,\triangle_{2}}$ for $\bsi_{h}$ in problem \eqref{DiscreV}, and piecewise linear space $P_{h,\triangle_{2}}$ for $u_{h}$. Triangulation $\mt$ is uniform. The numerical results are presented in Table \ref{Ta:L}. 
The convergence can still be observed in the L-shape domain. The converge rate is degenerate because the solution possesses singularities at the origin. Nevertheless, it is noted that the convergence rate of $\|u-u_{h}\|_{0}$ is higher than the other errors.

\begin{table}[h!]\footnotesize\tabcolsep 16pt
\begin{center}
\caption{The error and the order of convergence for the L-shape domain }\label{Ta:L}\vspace{-2mm}
\end{center}
\begin{center}
\begin{tabular}{c| c c c c c c }
\toprule  
      &$\|\sigma-\sigma_{h}\|_{0}$    &$h^{n}$     & $\|{\rd}{\bd}(\sigma-\sigma_{h})\|_{0}$   &$h^{n}$  &$\|u-u_{h}\|_{0}$  &$h^{n}$\\
\midrule  
1  & 3.0154e+00       &$-$        &1.1771e+02  &$-$        &2.7847e-01       &$-$       \\
2   &1.6652e+00      & 0.86      &4.9184e+01   &1.26     &4.8223e-02      & 2.53     \\
3   &1.1244e+00      &0.57       &2.1869e+01   & 1.17   & 2.1671e-02   &1.15          \\
4   &7.7274e-01       &0.54       &1.3617e+01  &0.68      &7.2229e-03      &1.59       \\
5   &5.3096e-01      &0.54        & 9.2420e+00   &0.56     &2.5820e-03       & 1.48    \\
\toprule
 &$\|Q_{h}u-u_{h}\|_{0}$    &$h^{n}$        &$|Q_{h}u-u_{h}|_{2,h,\triangle_{2}}$ &$h^{n}$  &$|u-u_{h}^{*}|_{2,h,\triangle_{2}}$  &$h^{n}$\\
\midrule
 1& 5.1725e-02  &$-$ &2.9152e-01  &$-$ &3.8066e+00 &$-$\\
2  &2.1221e-02    &1.29  &2.2187e-01  &0.39  &2.0648e+00 &0.88\\
  3 &9.5846e-03   &1.15 &1.6895e-01   &0.39 &1.3956e+00 &0.57\\
  4&4.4346e-03  &1.11 &1.1977e-01  &0.50 & 9.5956e-01 &0.54\\
 5 &2.0731e-03 & 1.10 &8.2704e-02    &0.53&6.5952e-01 &0.54\\
\bottomrule
\end{tabular}
\end{center}
\end{table}









\bibliographystyle{siam}
\bibliography{ref}

\begin{appendix}
\section{}
This appendix provides some ideas to construct the basis for $\Sigma_{k,\triangle_{d}}$. It is discussed for $k=3$ and $d=2$ while the ideas apply for $k\geq 3$ and $d=3$.

 For the case $d=2$ and $k=3$, let $\mathbf{x}_{1},\mathbf{x}_{2},\mathbf{x}_{3}$ be the vertices of a element $K\in \mt$. The affine mapping $F:\widehat{K}\rightarrow K$ reads
\[\mathbf{x}=F(\widehat{\mathbf{x}})=B\widehat{\mathbf{x}}+\mathbf{x}_{1},\]
with 
\[B=\begin{pmatrix}
\mathbf{x}_{2}-\mathbf{x}_{1},&\mathbf{x}_{3}-\mathbf{x}_{1}\\
\end{pmatrix}.\]
Suppose the triangle $\widehat{K}$ are spanned by $(0,0)$, $(1,0)$, $(0,1)$, and use $\widehat{\mathbf{x}}=(\widehat{x},\widehat{y})^{{\rm{T}}}$ for the vector in that coordinate.
Thus 
\begin{align}
\widehat{\mathbf{x}}_{1}=\begin{pmatrix}
0\\
0\\
\end{pmatrix}, ~~~~\widehat{\mathbf{x}}_{2}=\begin{pmatrix}
1\\
0\\
\end{pmatrix},~~~~\widehat{\mathbf{x}}_{3}=\begin{pmatrix}
0\\
1\\
\end{pmatrix}.
\end{align}

For each edge $e_{i}\in \mathcal{E}({K})$, the corresponding tangent vector is $\mathbf{t}_{i}=\mathbf{x}_{i-1}-\mathbf{x}_{i+1}$, $i=1,2,3$, where the indices are calculated ${\rm{mod} } \, 3$.  The unit outward normal vector of $e_{i}$ is denoted as $\mathbf{n}_{i}$. By the affine mapping,
\be
\label{af:nt}
\mathbf{n}_{i}=\frac{B^{-{\rm{T}}}\, \widehat{\mathbf{n}}_{i}}{|B^{-{\rm{T}}}\, \widehat{\mathbf{n}}_{i}|}, ~~~~\mathbf{t}_{i}=B\, \widehat{\mathbf{t}}_{i}.
\ee
The barycenter coordinates read
\begin{align}
\lambda_{2}&=\mathbf{n}_{2}\cdot(\mathbf{x}_{1}-\mathbf{x}),\\
\lambda_{3}&=\mathbf{n}_{3}\cdot(\mathbf{x}_{1}-\mathbf{x}),\\
\lambda_{1}&=1-\lambda_{2}-\lambda_{3}.
\end{align}
Denote ${\rm{J}}:={\rm{det}}(B)$. Note that ${\rm{J}}$ does not vanish at any point. Define for $\bta\in H({\rd}\bd,K;\ms)$, by the Piola transform \cite{2010Finite},

\be
\label{piola}
\bta(\mathbf{x}):=\frac{1}{{\rm{J}}}\, B\, \widehat{\bta}(\widehat{\mathbf{x}})\, B^{{\rm{T}}}.
\ee

Some fundamental properties of the Piola transform \eqref{piola} are presented in the subsequent lemmas.

\begin{lemma}
\label{refe}
If $\widehat{\bta}\in H({\bd},\widehat{K};\ms)$ satisfies $\widehat{\bta}\widehat{\mathbf{n}}|_{\partial \widehat{K}}=0$, then $\bta\in H({\bd},K;\ms)$ defined in \eqref{piola} satisfies $\bta\mathbf{n}|_{\partial K}=0$.
\end{lemma}
\begin{proof}
The combination of  \eqref{af:nt} and \eqref{piola} shows, on each edge $e\in \mathcal{E}(K)$, 
\be
\bta\mathbf{n}=\frac{1}{{\rm{J}}}B\widehat{\bta}B^{{\rm{T}}}\mathbf{n}=\frac{B\widehat{\bta}\widehat{\mathbf{n}}}{{\rm{J}}|B^{-{\rm{T}}}\widehat{\mathbf{n}}|}.
\ee
Thus $\widehat{\bta}\widehat{\mathbf{n}}|_{\widehat{e}}=0$ implies $\bta\mathbf{n}|_{e}=0$.
\end{proof}
\begin{lemma}
\label{refe1}
Suppose $\bta\in H({\rd}\bd,K;\ms)$, $q\in P_{k-1}(e)$, $e\in \mathcal{E}(K)$. If ${\rm{J}}>0$,  then 
\ben
(\bd\bta\cdot \mathbf{n}, q)_{e}=(\widehat{\bd}\widehat{\bta}\cdot \widehat{\mathbf{n}}, \widehat{q})_{\widehat{e}}.
\een
If ${\rm{J}}<0$, then
\ben
(\bd\bta\cdot \mathbf{n}, q)_{e}=-(\widehat{\bd}\widehat{\bta}\cdot \widehat{\mathbf{n}}, \widehat{q})_{\widehat{e}}.
\een
\end{lemma}
\begin{proof}
Since \cite{2010Finite} shows, $${\bd\bta}=\frac{1}{{\rm{J}}}B\, \widehat{\bd}\, \widehat{\bta}.$$ 
Let $q(x)=\widehat{q}(\widehat{x})$. It holds
\ben
(\bd\bta\cdot \mathbf{n},q)_{e}=(\frac{B\, \widehat{\bd}\, \widehat{\bta}B^{-{\rm{T}}}\widehat{\mathbf{n}}|e|}{{\rm{J}}|B^{-{\rm{T}}}\widehat{\mathbf{n}}||\widehat{e}|},\widehat{q})_{\widehat{e}}=(\frac{\widehat{\mathbf{n}}^{{\rm{T}}}\widehat{\bd}\, \widehat{\bta}|e|}{{\rm{J}}|B^{-{\rm{T}}}\widehat{\mathbf{n}}||\widehat{e}|},\widehat{q})_{\widehat{e}}=\frac{|{\rm{J}}|}{{\rm{J}}}(\widehat{\bd}\widehat{\bta}\cdot \widehat{\mathbf{n}},\widehat{q})_{\widehat{e}}.
\een
This concludes the proof.
\end{proof}

The basis for $\Sigma_{k,\triangle_{2}}$ are formed as follows. For $k=3$, only the degrees of freedom \eqref{2dDof:1}--\eqref{2dDof:3} are adopted. The first step is to construct basis functions for the degrees of freedom \eqref{2dDof:3}, which are denoted by $\bta_{h,i}$, $i=1, 2, \cdots, 9$.

Recall the $H(\bd,K;\ms)$ bubble functions $\bm{\vartheta}_{h}$ introduced in \cite{Hu2014A},
\ben
\bm{\vartheta}_{h}\in\Sigma_{K,b}:=\sum_{1\leq i\leq 3}\lambda_{i-1}\lambda_{i+1}P_{1}(K)\mathbf{t}_{i}\mathbf{t}_{i}^{{\rm{T}}}
\een
with
\ben
\bm{\vartheta}_{h}\mathbf{n}_{j}|_{e_{j}}=\mathbf{0}, ~~~~j=1,2,3.
\een

Lemma \ref{refe} ensures that $\bta_{h,i}$ can be obtained from the basis functions $\widehat{\bta}_{h,i}$ defined on the reference element $\widehat{K}$.  Let the nine basis functions of $\widehat{\Sigma}_{\widehat{K},b}$ be $\bm{\widehat{\vartheta}}_{h,i}$, $i=1,2,\cdots,9$. To be precise, 
\ben
\begin{aligned}
\bm{\widehat{\vartheta}}_{h,1}&=\frac{9}{2}\widehat{\lambda}_{2}\widehat{\lambda}_{3}(3\widehat{\lambda}_{2}-1)\widehat{\mathbf{t}}_{1}\widehat{\mathbf{t}}_{1}^{{\rm{T}}};\\
\bm{\widehat{\vartheta}}_{h,2}&=\frac{9}{2}\widehat{\lambda}_{3}\widehat{\lambda}_{2}(3\widehat{\lambda}_{3}-1)\widehat{\mathbf{t}}_{1}\widehat{\mathbf{t}}_{1}^{{\rm{T}}};\\
\bm{\widehat{\vartheta}}_{h,3}&=\frac{9}{2}\widehat{\lambda}_{3}\widehat{\lambda}_{1}(3\widehat{\lambda}_{3}-1)\widehat{\mathbf{t}}_{2}\widehat{\mathbf{t}}_{2}^{{\rm{T}}};\\
\bm{\widehat{\vartheta}}_{h,4}&=\frac{9}{2}\widehat{\lambda}_{1}\widehat{\lambda}_{3}(3\widehat{\lambda}_{1}-1)\widehat{\mathbf{t}}_{2}\widehat{\mathbf{t}}_{2}^{{\rm{T}}};\\
\bm{\widehat{\vartheta}}_{h,5}&=\frac{9}{2}\widehat{\lambda}_{1}\widehat{\lambda}_{2}(3\widehat{\lambda}_{1}-1)\widehat{\mathbf{t}}_{3}\widehat{\mathbf{t}}_{3}^{{\rm{T}}};\\
\bm{\widehat{\vartheta}}_{h,6}&=\frac{9}{2}\widehat{\lambda}_{2}\widehat{\lambda}_{1}(3\widehat{\lambda}_{2}-1)\widehat{\mathbf{t}}_{3}\widehat{\mathbf{t}}_{3}^{{\rm{T}}};\\
\bm{\widehat{\vartheta}}_{h,7}&=27\widehat{\lambda}_{1}\widehat{\lambda}_{2}\widehat{\lambda}_{3}\widehat{\mathbf{t}}_{1}\widehat{\mathbf{t}}_{1}^{{\rm{T}}};\\
\bm{\widehat{\vartheta}}_{h,8}&=27\widehat{\lambda}_{1}\widehat{\lambda}_{2}\widehat{\lambda}_{3}\widehat{\mathbf{t}}_{2}\widehat{\mathbf{t}}_{2}^{{\rm{T}}};\\
\bm{\widehat{\vartheta}}_{h,9}&=27\widehat{\lambda}_{1}\widehat{\lambda}_{2}\widehat{\lambda}_{3}\widehat{\mathbf{t}}_{3}\widehat{\mathbf{t}}_{3}^{{\rm{T}}}.
\end{aligned}
\een
Assume 
\ben
\widehat{\bta}_{h,i}=\sum_{i=1}^{9}\alpha_{j}^{(i)}\bm{\widehat{\vartheta}}_{h,j},~~~~ \alpha_{j}^{i}\in\mathbb{R}.
\een

The corresponding basis functions for degrees of freedom \eqref{2dDof:3} on $\widehat{K}$ can be calculated immediately. Suppose $C$ denotes the $9\times 9$ coefficients matrix consisting of $\alpha_{j}^{(i)}$, and let $C(i, j)=\alpha_{j}^{(i)}$, then 
\ben
C=\begin{pmatrix}
       0              &0             & 4/9            &2/9          &  2/3            &4/3            &1/3          & -1/3           & 1/9     \\
       0              &0             & 4/3           & 2/3          &  2/9           & 4/9            &1/9          & -1/3            &1/3     \\
       0              &0            &-40/9        &  -20/9        &  -20/9         & -40/9         & -10/9           &20/9        &  -10/9    \\ 
      -4/3          & -8/3        &    0           &   0            & -4/9           &-2/9            &2/9            &0            & -1/3     \\
      -4/9           &-8/9         &   0           &   0            & -4/3           &-2/3            &2/9          & -2/9          & -1/9    \\ 
      40/9          & 80/9        &    0           &   0            & 40/9        &   20/9         & -20/9          & 10/9        &   10/9   \\  
      -8/9          & -4/9          & -2/3         &  -4/3         &   0             & 0             &-1/9          & -2/9          &  2/9   \\  
      -8/3           &-4/3          & -2/9          & -4/9         &   0             & 0             &-1/3           & 0             & 2/9     \\
      80/9           &40/9        &   20/9        &   40/9      &      0            &  0           &  10/9          & 10/9          &-20/9   \\     
\end{pmatrix}.
\een
This leads to 
\begin{align*}
\widehat{\bta}_{h,1}=\begin{pmatrix}
9\widehat{\lambda}_{1}\widehat{\lambda}_{2}^{2} &-9\widehat{\lambda}_{1}\widehat{\lambda}_{2}\widehat{\lambda}_{3}\\
-9\widehat{\lambda}_{1}\widehat{\lambda}_{2}\widehat{\lambda}_{3} &3\widehat{\lambda}_{1}\widehat{\lambda}_{3}^{2}\\
\end{pmatrix}, ~\widehat{\bta}_{h,2}=\begin{pmatrix}
3\widehat{\lambda}_{1}\widehat{\lambda}_{2}^{2} &-9\widehat{\lambda}_{1}\widehat{\lambda}_{2}\widehat{\lambda}_{3}\\
-9\widehat{\lambda}_{1}\widehat{\lambda}_{2}\widehat{\lambda}_{3} &9\widehat{\lambda}_{1}\widehat{\lambda}_{3}^{2}\\
\end{pmatrix},
\end{align*}
\begin{align*}
\widehat{\bta}_{h,3}=\begin{pmatrix}
-30\widehat{\lambda}_{1}\widehat{\lambda}_{2}^{2} &60\widehat{\lambda}_{1}\widehat{\lambda}_{2}\widehat{\lambda}_{3}\\
60\widehat{\lambda}_{1}\widehat{\lambda}_{2}\widehat{\lambda}_{3} &-30\widehat{\lambda}_{1}\widehat{\lambda}_{3}^{2}\\
\end{pmatrix},
\end{align*}
\begin{align*}
\widehat{\bta}_{h,4}=\begin{pmatrix}
-3\widehat{\lambda}_{2}(\widehat{\lambda}_{2}^{2}+8\widehat{\lambda}_{2}\widehat{\lambda}_{3}-\widehat{\lambda}_{2}+10\widehat{\lambda}_{3}^{2}-7\widehat{\lambda}_{3}+\widehat{\lambda}_{1}) &9\widehat{\lambda}_{2}\widehat{\lambda}_{3}(\widehat{\lambda}_{3}-\widehat{\lambda}_{1})\\
9\widehat{\lambda}_{2}\widehat{\lambda}_{3}(\widehat{\lambda}_{3}-\widehat{\lambda}_{1}) &-9\widehat{\lambda}_{2}\widehat{\lambda}_{3}^{2}\\
\end{pmatrix},
\end{align*}
\begin{align*}
\widehat{\bta}_{h,5}=\begin{pmatrix}
-3\widehat{\lambda}_{2}(3\widehat{\lambda}_{2}^{2}+12\widehat{\lambda}_{2}\widehat{\lambda}_{3}-3\widehat{\lambda}_{2}+10\widehat{\lambda}_{3}^{2}-9\widehat{\lambda}_{3}+3\widehat{\lambda}_{1}) &3\widehat{\lambda}_{2}\widehat{\lambda}_{3}(\widehat{\lambda}_{3}-3\widehat{\lambda}_{1})\\
3\widehat{\lambda}_{2}\widehat{\lambda}_{3}(\widehat{\lambda}_{3}-3\widehat{\lambda}_{1}) &-3\widehat{\lambda}_{2}\widehat{\lambda}_{3}^{2}\\
\end{pmatrix},
\end{align*}
\begin{align*}
\widehat{\bta}_{h,6}=\begin{pmatrix}
30\widehat{\lambda}_{2}(\widehat{\lambda}_{2}^{2}+6\widehat{\lambda}_{2}\widehat{\lambda}_{3}-\widehat{\lambda}_{2}+6\widehat{\lambda}_{3}^{2}-5\widehat{\lambda}_{3}+\widehat{\lambda}_{1}) &30\widehat{\lambda}_{2}\widehat{\lambda}_{3}(2\widehat{\lambda}_{1}-\widehat{\lambda}_{3})\\
30\widehat{\lambda}_{2}\widehat{\lambda}_{3}(2\widehat{\lambda}_{1}-\widehat{\lambda}_{3}) &30\widehat{\lambda}_{1}\widehat{\lambda}_{2}^{2}\\
\end{pmatrix},
\end{align*}
\begin{align*}
\widehat{\bta}_{h,7}=\begin{pmatrix}
-3\widehat{\lambda}_{2}^{2}\widehat{\lambda}_{3} & 3\widehat{\lambda}_{2}\widehat{\lambda}_{3}(2\widehat{\lambda}_{2}-\widehat{\lambda}_{1})\\
3\widehat{\lambda}_{2}\widehat{\lambda}_{3}(2\widehat{\lambda}_{2}-\widehat{\lambda}_{1}) &-3\widehat{\lambda}_{3}(10\widehat{\lambda}_{2}^{2}+12\widehat{\lambda}_{2}\widehat{\lambda}_{3}-9\widehat{\lambda}_{2}+3\widehat{\lambda}_{3}^{2}-3\widehat{\lambda}_{3}+3\widehat{\lambda}_{1})\\
\end{pmatrix},
\end{align*}
\begin{align*}
\widehat{\bta}_{h,8}=\begin{pmatrix}
-9\widehat{\lambda}_{2}^{2}\widehat{\lambda}_{3} & 9\widehat{\lambda}_{2}\widehat{\lambda}_{3}(2\widehat{\lambda}_{2}-\widehat{\lambda}_{1})\\
9\widehat{\lambda}_{2}\widehat{\lambda}_{3}(2\widehat{\lambda}_{2}-\widehat{\lambda}_{1}) &-3\widehat{\lambda}_{3}(10\widehat{\lambda}_{2}^{2}+8\widehat{\lambda}_{2}\widehat{\lambda}_{3}-7\widehat{\lambda}_{2}+\widehat{\lambda}_{3}^{2}-\widehat{\lambda}_{3}+\widehat{\lambda}_{1})\\
\end{pmatrix},
\end{align*}
\begin{align*}
\widehat{\bta}_{h,9}=\begin{pmatrix}
30\widehat{\lambda}_{2}^{2}\widehat{\lambda}_{3} & 30\widehat{\lambda}_{2}\widehat{\lambda}_{3}(2\widehat{\lambda}_{1}-\widehat{\lambda}_{2})\\
30\widehat{\lambda}_{2}\widehat{\lambda}_{3}(2\widehat{\lambda}_{1}-\widehat{\lambda}_{2}) &30\widehat{\lambda}_{3}(6\widehat{\lambda}_{2}^{2}+6\widehat{\lambda}_{2}\widehat{\lambda}_{3}-5\widehat{\lambda}_{2}+\widehat{\lambda}_{3}^{2}-\widehat{\lambda}_{3}+\widehat{\lambda}_{1})\\
\end{pmatrix}.
\end{align*}
Hence $\bta_{h,i}$, $i=1,2,\cdots,9$ follow by Piola transform \eqref{piola}. These basis $\bta_{h,i}$, $i=1,2,\cdots,9$, satisfy 
\ben 
\begin{aligned}
\bta_{h,i}\mathbf{n}_{j}|_{e_{j}}&=0~~~~&1\leq i\leq 9,~j=1,2,3,\\
d_{i,e}(\bta_{h,j})&=\delta_{ij}~~~~&1\leq i,j \leq 9.\\
\end{aligned}
\een
Here $d_{i,e}(\cdot)$, $i=1,2,\cdots,9$ are defined by
\ben 
\begin{aligned}
d_{1,e}(\cdot)=\int_{e_{1}}\mathbf{n}_{1}^{{\rm{T}}}\bd(\cdot)\lambda_{2},~~d_{2,e}(\cdot)=\int_{e_{1}}\mathbf{n}_{1}^{{\rm{T}}}\bd(\cdot)\lambda_{3},~~d_{3,e}(\cdot)=\int_{e_{1}}\mathbf{n}_{1}^{{\rm{T}}}\bd(\cdot)\lambda_{2}\lambda_{3},\\
d_{4,e}(\cdot)=\int_{e_{2}}\mathbf{n}_{2}^{{\rm{T}}}\bd(\cdot)\lambda_{3},~~d_{5,e}(\cdot)=\int_{e_{2}}\mathbf{n}_{2}^{{\rm{T}}}\bd(\cdot)\lambda_{1},~~d_{6,e}(\cdot)=\int_{e_{2}}\mathbf{n}_{2}^{{\rm{T}}}\bd(\cdot)\lambda_{3}\lambda_{1},\\
d_{7,e}(\cdot)=\int_{e_{3}}\mathbf{n}_{3}^{{\rm{T}}}\bd(\cdot)\lambda_{1},~~d_{8,e}(\cdot)=\int_{e_{3}}\mathbf{n}_{3}^{{\rm{T}}}\bd(\cdot)\lambda_{2},~~d_{9,e}(\cdot)=\int_{e_{3}}\mathbf{n}_{3}^{{\rm{T}}}\bd(\cdot)\lambda_{1}\lambda_{2}.\\
\end{aligned}
\een

The second step is to construct the remainder $21$ basis functions $\bta_{h,i}$, $i=10,11,\cdots,30$, for $\Sigma_{3,\triangle_{2}}$. These basis satisfy
\ben
\int_{e}\mathbf{n}_{e}^{{\rm{T}}}\bd\bta_{h,i}\, p_{2}\, ds=0~~ \text{for all}~~p_{2}\in P_{2}(e), \,  e\in \mathcal{E}(K).
\een
Similarly, recall the rest two types of basis functions in \cite{Hu2014A}, which are vertex-based basis functions and edge-based basis functions with nonzero fluxes. On element $K\in \mt$, the remainder $21$ basis functions of $\Sigma_{3,\triangle_{2}}$ can be derived from the following two classes of basis functions in \cite{Hu2014A}. 
\begin{enumerate}
\item Vertex-based basis functions. The $9$ basis functions  in \cite{Hu2014A} are defined by
\begin{align*}
\bm{{\varphi}}_{h,i}=\phi_{1}\mathbb{T}_{i},~~~~i=1,2,3,\\
\bm{{\varphi}}_{h,i+3}=\phi_{2}\mathbb{T}_{i},~~~~i=1,2,3,\\
\bm{{\varphi}}_{h,i+6}=\phi_{3}\mathbb{T}_{i},~~~~i=1,2,3,
\end{align*}
with the Lagrange nodal basis functions in $P_{3}(K)$
\[\phi_{i}=\frac{1}{2}\lambda_{i}(3\lambda_{i}-1)(3\lambda_{i}-2),~~~~i=1,2,3,\]
and
\ben
\mathbb{T}_{1}:=\begin{pmatrix}
1 &0\\
0 &0\\
\end{pmatrix},~~~~\mathbb{T}_{2}:=\begin{pmatrix}
0 &1\\
1&0\\
\end{pmatrix},~~~~\mathbb{T}_{3}:=\begin{pmatrix}
0 &0\\
0&1\\
\end{pmatrix}.
\een
\item Edge-based basis functions with nonzero fluxes. For $e_{i}\in \mathcal{E}(K)$, $i=1,2,3$, denote
$$\mathbf{x}_{e_{i},1}=\frac{1}{3}(2\mathbf{x}_{i+1}+\mathbf{x}_{i-1}),~~\mathbf{x}_{e_{i},2}=\frac{1}{3}(2\mathbf{x}_{i-1}+\mathbf{x}_{i+1}).$$ 
The associated Lagrange nodal basis functions are
\ben
\phi_{e_{i},1}=\frac{9}{2}\lambda_{i+1}\lambda_{i-1}(3\lambda_{i+1}-1),~~ \phi_{e_{i},2}=\frac{9}{2}\lambda_{i-1}\lambda_{i+1}(3\lambda_{i-1}-1).
\een
The $12$ edge-based basis functions $\bm{{\varphi}}_{h,10},\cdots,\bm{{\varphi}}_{h,21}$ (with nonzero fluxes) in \cite{Hu2014A} are 
\begin{align*}
\phi_{e_{i},j}\mathbf{n}_{i}\mathbf{n}_{i}^{{\rm{T}}}, ~~\frac{1}{2}\phi_{e_{i},j}(\mathbf{t}_{i}\mathbf{n}_{i}^{{\rm{T}}}+\mathbf{n}_{i}\mathbf{t}_{i}^{{\rm{T}}}), ~~~~i=1,2,3,~j=1,2,
\end{align*}
respectively.
\end{enumerate}

The basis functions of $\Sigma_{h,\triangle_{2}}$ have the forms
\ben
\bta_{h,i+9}=\bm{{\varphi}}_{h,i}-\sum_{j=1}^{9}\beta_{j}^{(i)}{\bta}_{h,j},~~~~ i=1,2,\cdots,21.
\een
The coefficients $\beta_{1}^{(i)},\cdots,\beta_{9}^{(i)}$ are constants, given by 
\begin{align*}
(\mathbf{n}_{l}^{{\rm{T}}}\bd \bm{{\varphi}}_{h,i},\lambda_{l+1})_{e_{l}},~ (\mathbf{n}_{l}^{{\rm{T}}}\bd \bm{{\varphi}}_{h,i},\lambda_{l-1})_{e_{l}},~(\mathbf{n}_{l}^{{\rm{T}}}\bd \bm{{\varphi}}_{h,i},\lambda_{l-1}\lambda_{l+1})_{e_{l}}, ~~l=1,2,3,
\end{align*}
respectively.

\begin{remark}
The implementation of $\bta_{h,i}$, $i=10,11,\cdots,30$  can also rely on the reference element $\widehat{K}$. For instance, according to Lemma \ref{refe1}, 
\be
(\mathbf{n}_{1}^{{\rm{T}}}\bd \bm{{\varphi}}_{h,1},\lambda_{2})_{e_{l}}=(\widehat{\mathbf{n}}_{1}^{{\rm{T}}}\widehat{\bd} \widehat{\bm{{\varphi}}}_{h,1},\widehat{\lambda}_{2})_{\widehat{e}_{1}},
\ee
where  
\ben
\widehat{\bm{{\varphi}}}_{h,1}={\rm{J}}\, B^{-1}\bm{{\varphi}}_{h,1}B^{-{\rm{T}}}={\rm{J}}\, B^{-1}\phi_{1}\mathbb{T}_{1}B^{-{\rm{T}}}:=\phi_{1}M_{1}.
\een
Note that matrix $M_{1}={\rm{J}}\, B^{-1}\mathbb{T}_{1}B^{-{\rm{T}}}$. This shows
\ben
\begin{split}
(\mathbf{n}_{1}^{{\rm{T}}}\bd \bm{{\varphi}}_{h,1},\lambda_{2})_{e_{1}}&=(\widehat{\mathbf{n}}_{1}^{{\rm{T}}}\widehat{\bd}(\phi_{1}M_{1}),\widehat{\lambda}_{2})_{\widehat{e}_{1}}\\
&=(\widehat{\mathbf{n}}_{1}^{{\rm{T}}}M_{1}\widehat{\nabla}\phi_{1},\widehat{\lambda}_{2})_{\widehat{e}_{1}}\\
&=\widehat{\mathbf{n}}_{1}^{{\rm{T}}}M_{1}(\widehat{\nabla}\phi_{1},\widehat{\lambda}_{2})_{\widehat{e}_{1}}.
\end{split}
\een
The term $(\widehat{\nabla}\phi_{1},\widehat{\lambda}_{2})_{\widehat{e}_{1}}$ can be calculated exactly on $\widehat{K}$. Then some transformations lead to $ (\mathbf{n}_{1}^{{\rm{T}}}\bd \bm{{\varphi}}_{h,1},\lambda_{2})_{e_{1}}$ which is $\beta_{1}^{1}$.
\end{remark}

\end{appendix}

\end{document}